\documentclass[a4paper]{article}

\newtheorem{thm}{Theorem}[section]
\newtheorem{lem}[thm]{Lemma}
\newtheorem{pro}[thm]{Proposition}
\newtheorem{cor}[thm]{Corollary}
\newtheorem{dfn}[thm]{Definition}
\begin{document} 
\setlength{\parindent}{.4in}
\setlength{\parskip}{2.2mm}

\setcounter{section}{0}

\title{Projective tensor products and $A_p^q$ spaces}
\author{William Moran and H.Kumudini Dharmadasa}
\maketitle

\section*{Introduction}

Let $G$ be a group, let $H$ and $K$ be two subgroups
 of $G,$ and let ${\pi }$ and ${\gamma }$ be representations of $H$ and $K,$
respectively. If $G$ is finite, Mackey's results assert that the \sf intertwining 
number \rm  of the two \sf induced representations \rm $U^{\pi }$ and
 $U^{\gamma } $ of $G$   can be expressed as a sum of intertwining
 numbers of
 the representations ${\pi ^x}$ and ${\gamma^y }$ of the subgroups 
$H^x \cap K^y $, $x,y \in G.$ In the case of an infinite group, if the 
subgroups are open and closed, a similar characterization 
is possible especially when  ${\pi }$ and ${\gamma }$ are one-dimensional.
If the subgroups are closed, Mackey showed that the above criteria for
computing the intertwining number holds for the space of those operators 
 which are in 
the Hilbert-Schmidt class.

%The fact that every continuous linear map of an $L_1 $ space into a separable reflexive space 
%can be better represented as an integral operator led Moore \cite{moo} to extend the 
%concept of induced representations to include the action of a group on a Banach 
%space by isometries. He proved that the Frobenius
%Reciprocity Theorem remains true under these modifications and the
% assumption that the 
%corresponding $G$-coset space possesses an invariant measure. 

Among other developments that are important for us, the first is the work of 
Rieffel\cite{rief2} on Banach G-modules and their products. He proved,
 in particular, that 
 $$(V\otimes_S W)^* \cong Hom_S (V,W^*),$$

\noindent where $S$ is a set, $V$ and $W$ are two $S$-modules,
$\otimes_S $ denotes the projective tensor product of $V$ and $W$
and $Hom_S (V,W^*)$ is  the space of intertwining operators of 
the Banach G-modules. Applying this to $L_p (G) $ spaces 
($1\leq p \leq \infty $)
of complex-valued functions defined on a 
group $G$, Rieffel obtained the result that, under certain conditions, the 
corresponding intertwining operators (multipliers) form the dual space of 
the space of functions $A_p^q$: a subset of an $L_r $ space
 (where $r$ is related to $p$ and $q$ as described in  Prop.3.18) consisting 
of those functions which can be written as a sum of convolution of functions
from $L_p$ and $L_q.$ This is the context in which we 
shall set our study of intertwining operators, that is, regarding the space
 of such operators as the dual of a tensor product space.    

 Herz
\cite{her}
 studied the predual of the space of intertwining operators of the regular 
representations of $G$ on $L_p$ and $L_q$ and was able to show, in 
particular,
that the tensor product space is an algebra of functions on $G$ and, in some 
sense, a natural analogue of the space of absolutely convergent Fourier Series.
Our aim is to extend the Herz- Rieffel results from regular representations
which may be seen as induced representation from the trivial subgroup to
arbitrary induced representations.

In order to complete this analysis we shall need to go beyond spaces of functions
on $G$ to sections of Banach (semi-)bundles on $G$. The concept of a Banach 
bundle was developed by Fell in 1977 and we shall use it as the 
appropriate device for the study of the tensor product spaces. Unfortunately,
in the most general case, our semi-bundle will fail to be a bundle in the 
complete sense, but will be more akin to the objects studied by Dauns and 
Hofmann\cite{dau}.

\section {Preliminaries}
\subsection{$\lambda$-functions}
We shall assume throughout that all the topological spaces under 
consideration are second countable.  

Let $G$ be a locally compact topological group. We denote the right-invariant 
Haar measure on $G$ by $\nu_G.$ $e$ denotes the identity element of the 
group. For a subgroup $H$ of $G,$ the canonical mapping from $G$ to
 the set of right-cosets $G/H$ is denoted by $p_H .$

 A  real-valued 
 function $\rho_H $  defined on $G$ which satisfies
\begin{eqnarray}
\rho_H (hx)&= &({\Delta_H (h)/ {\Delta_G (h)} }) \rho_H (x),
\end{eqnarray}
\noindent where  $x\in G$ and $h\in H,$
is called a $\rho$-function.
 The existence of a strictly positive continuous $\rho$-function
has been established in a number of places in the literature 
(see Mackey\cite{mack1}, Gaal\cite{gal}). In particular,
it is known that for every closed subgroup $H$ in $G$ there exists a function 
$\beta$ on $G$ with $\int_H \beta (hx) d\nu_H (h) =1$ for all $x\in G $
which gives rise to a $\rho$-function of the required nature. The details 
of such a $\beta$ function are given in the following Lemma.

\begin{lem} 
 For every closed subgroup $H$ of a locally compact group 
$G,$ there exists a function $\beta $ on $G$ with the following properties:
\begin{description}
\item[(a)] if $K$ is any compact set in $G,$ then $\beta $ coincides on the strip 
$HK$ with a function in $C_0^{+} (G);$

\item[(b)] $\int_H \beta (hx) d\nu_H (h) = 1 $ for all $x\in G.$
\end{description}
\end{lem}

\noindent Proof: See Reiter\cite{rei}, Chapter 8, section 1.9. 
\begin{flushright} \large$\diamond$\\
\end{flushright}
 %\large$\diamond$\\

 A function 
 $\beta $ on $G$ satisfying the properties stated in Lemma 1.1 
is called a \sf Bruhat function \rm for $H.$

 Given a Bruhat function $\beta  $ for a closed subgroup $H,$ a 
$\rho$-function can be obtained by letting
$$\rho_H (x) = \int_H \beta (hx) \Delta_G (h)\Delta_H (h^{-1} ) d\nu_H (h) .$$
Then $\rho_H $ is continuous (cf. (a) and  \cite{rei}, Chapter 3, section 3.2,
Remark) and strictly positive for all $x\in G.$ 
 
 For a given $\rho$-function ${ \rho (sy)/ { \rho (s)}}$
is a Borel function of $s$ and $y$ which is constant on the right $H\times G$
cosets in $G\times G.$ Since there is a natural homeomorphism from this
coset space to $(G/H) \times G,$
these $\rho$-functions give rise to a
unique Borel function $\lambda_{\rho} $ on $(G/H)\times G$ such that 
$$\lambda_{\rho}  (p_H (s), y) ={ \rho (sy)\over { \rho (s)}}$$ for all $s$ and
$y$ in $G.$ 
This function $\lambda_{\rho} $ has the following properties:
\begin{description}
\item[(a)] for all $x\in (G/H)$ and $s,t\in G,\quad \lambda_{\rho}(x,st)=
\lambda_{\rho}(x.s,t) \lambda_{\rho} (x,s);$
\item[(b)] for all $h\in H, \lambda_{\rho}(p_H (e),h) =  {\Delta_H (h)/ {\Delta_G
 (h)} };$ 
\item[(c)] $\lambda_{\rho}(p_H(e),t) $ is bounded on compact sets as a 
function of $t.$
\end{description}
(See, for example, Gaal\cite{gal}, p.263, Lemma 10.) For a given measure 
$\mu $ on $G/H$ and $y\in G/H,$ let $\mu_y$ denote the translated measure on 
$G/H$ defined by $\mu_y (E) = \mu ([E]y).$ 
 It is well known that for a given arbitrary $\rho$-function on $G$
 there  exists a quasi-invariant measure $\mu $ in the right coset space
$G/H$ such that for all
$y\in G,$ the corresponding   $\lambda$-function $\lambda_{\rho }$ has the
property that $\lambda_{\rho } (\cdot ,y)$ is a Radon-Nikodym derivative of the measure
$\mu_y$ with respect to the measure $\mu.$  Any two
have the  same null sets and hence are mutually absolutely continuous. A Borel
set $E$ in $G/H$ is a null
 set if and only if $p_H^{-1}(E) $ has Haar measure zero.  Let us write  $\mu \succ \lambda $ to mean  that
for all $y\in G,$ $\lambda (\cdot ,y) $ is a Radon-Nikodym derivative of the
measure  $\mu_y$ with respect to $\mu .$ 
 The relations $\mu \succ
\lambda $ and $ \lambda =  \lambda_{\rho}$ between quasi-invariant measures,
$\lambda$-functions and $\rho $-functions have the following properties:
\begin{description}
\item[(i)]  Every $ \lambda $-function is of the form $ \lambda_{\rho}$; $ \lambda_{\rho_1}=
 \lambda_{\rho_2}$ if and only if ${{\rho_1}/ {\rho_2}}$ is a constant.

\item[(ii)] If $\mu_1 \succ \lambda $  and $\mu_2 \succ \lambda $ 
then $\mu_1 $ is a constant multiple of $\mu_2 .$

\item[(iii)]  If $\mu \succ \lambda_1 $ and $\mu \succ \lambda_2 $ 
then for all $t,$ $\lambda_1 (\cdot ,t)$ = $\lambda_2 (\cdot ,t)$ almost everywhere 
in $G/H.$
\end{description}
( See Mackey\cite{mack1} for a detailed
study on ${\rho }$-functions and related  $\lambda $-functions).

 The quasi-invariant measure on the homogeneous space $G/H$ of a subgroup 
$H$ of
a group $G$ will be denoted by   $\mu_{H}$ and the Radon-Nikodym derivative of the 
measure $E\mapsto \mu_{H} ([E]y)$ with respect to the measure $\mu_{H}$
is denoted by $\lambda_{H}(. ,y). $

\sf For  simplicity of notation, $\lambda_H (p_H (x), y)$ will be
 written as $\lambda_H (x,y),$ or by $\lambda (x,y)$ if the subgroup $H$
 is clearly understood.\rm 

The following result, which appears in several places in the literature, 
is of fundemental importance for our purposes.
\begin{cor} For $x\in G$ let $\stackrel{.}{x}= p_H (x).$ If $\mu $ denotes 
the quasi-invariant measure corresponding to the function $\rho $ then
$$\int_G f(x)\rho (x) d\nu_G (x) = \int_{G\over H} \int_{H} f(hx) d\nu_H (h) 
d\mu (\stackrel{.}{x}), \quad f\in C_0 (G).$$
\end{cor}

\noindent Proof: See, for example, Gaal\cite{gal}, p.263, Corollary to Theorem 9.

\begin{flushright} \large$\diamond$\\
\end{flushright}

 \subsection {Banach Bundles}
 Here we recall definitions and a few results in terms of Banach  bundles.
(see \cite{fel1}, Chapter 2 and \cite{fel2} for further details.) 

  A \sf bundle \rm ${\underline {\cal B}} $ over a 
Hausdorff space $X$ is a pair $ ({\cal B}, \theta)$
 such that ${\cal B}$ is a Hausdorff space called the 
\sf bundle space \rm of ${\underline {\cal B}}$ and $ \theta : 
{\cal B} \mapsto X$ is a continuous open
surjection called the \sf bundle projection \rm of ${\underline {\cal B}}.$
 $X$ is called the \sf base  space \rm
of ${\underline {\cal B}},$ and for $x\in X ,$  $\theta ^{-1} (x) = \{
\xi : \theta (\xi ) = x, \xi\in {\cal B}\}$  is 
called the \sf fibre \rm over $X$ and is denoted by ${\cal B}_x . $

A bundle ${\underline {\cal B}}=({\cal B}, \theta)$ over $X$ is a 
\sf Banach semi-bundle  \rm over $X$ if we can define a 
norm making each fibre ${\cal B}_x $  into a Banach space satisfying the 
following conditions:
\begin{description}
\item[(a)] $\xi \mapsto \| \xi \|$ is upper semi-continuous 
 on ${\cal B}$ to ${\cal R}.$
\item[(b)] The operation $+$ is continuous on the set 
$\{(\xi ,\eta )\in {\cal B}\times 
{\cal B} : \theta (\xi ) 
=\theta (\eta )\}$ to ${\cal B}.$
\item[(c)] For each $\lambda $ in ${{\cal C}},$ the map $\xi \mapsto
 \lambda  .\xi $ 
is continuous on ${\cal B}$ to ${\cal B}.$
\item[(d)] If $x\in X$ and $\{\xi_i \}$ is a net of elements of ${\cal B}$ such
 that $\|\xi_i \|
\rightarrow 0$ and $\theta (\xi_i ) \rightarrow x,$ then $\xi_i \rightarrow
 0_x,$ where $0_x$
denotes the zero element of the Banach space  ${\cal B}_x .$
\end{description}

A bundle ${\underline {\cal B}}=({\cal B}, \theta)$ over $X$ is called a 
\sf Banach  bundle \rm if it
 satisfies (b), (c) and (d) above together with the condition that
\begin{description}
\item[(\~a)] $\xi \mapsto \| \xi \|$ is continuous
 on ${\cal B}$ to ${\cal R}.$
\end{description}

Given a Banach space $A$ and a Hausdorff space $X,$ it is easy to construct
a Banach bundle by letting ${\cal B}= A\times X$ and $\theta (\xi, x)= x .$
Then $({\cal B}, \theta )$ is a bundle over $X$ and if we equip each
 fibre $A\times \{ x\} $ with the Banach space structure making $\xi \mapsto 
(\xi , x)$ an isometric isomorphism, then it becomes a Banach bundle.  
The Banach bundle $({\cal B}, \theta )$ so constructed is called a \sf
trivial Banach bundle. \rm 

 Let $X$ and $Y$ be any two Hausdorff spaces and $\phi : Y \mapsto X$ 
be a continuous map.
Suppose  ${\underline {\cal B}}=({\cal B}, \theta)$ is a Banach (semi-)bundle 
over $X.$ Let ${\cal B}^{\# }$ be the topological 
subspace $\{ (y,\xi ): y\in Y , \xi \in {\cal B}, \phi (y) = \theta (\xi )\}
$ of $Y\times {\cal B}$; and 
define $\theta^{\# } : {\cal B}^{\# } \mapsto Y $ by $\theta^{\# }
 (y,\xi ) = y.$ Then $\theta^{\# } $ is a continuous open surjection 
since $\theta $ is open. Hence $({\cal B}^{\# }, \theta^{\# } )$ is a bundle over
 $Y.$ For $y\in Y,$ we make 
${\cal B}^{\# }_y = {\theta^{\# }}^{-1}(y)$ into a Banach space in such a way that the bijection 
$\xi \mapsto (y,\xi )$
of ${\cal B}_{\phi (y)}$ onto ${\cal B}^{\# }_y$ becomes a linear isometry.
 Then $({\cal B}^{\# },\theta^{\# } ), $ denoted by 
${\underline {\cal B}}^{\# } ,$ becomes a Banach (semi-)bundle
 which is called the \sf Banach (semi-)bundle retraction \rm of
 ${\underline {\cal B}}$ by $\phi .$

 Let $i^{\# }:{\cal B}^{\# }\mapsto {\cal B}$ be the surjection 
 given by $i^{\# } (y,\xi )= \xi.$ Then, we have the following diagram:

$${\cal B}^{\# }\stackrel{{i^{\# }}}{\longmapsto }{\cal B}$$
$$\theta^{\# } \downarrow  \quad \quad \quad \quad \downarrow \theta $$
$$Y\stackrel{\phi}{\mapsto  } X$$
\noindent Since
$\theta (i^{\# } (y,\xi ))=\theta (\xi ) =\phi (y) = \phi (\theta^{\# } 
(y, \xi )),$ 
we have $\theta i^{\# } = \phi \theta^{\# }$, and the diagram commutes.

Suppose ${\underline {\cal B}}=({\cal B}, \theta)$ and ${\underline {\cal D}
}=({\cal D}, \vartheta)$ are Banach (semi-)bundles over the same base space
 $X.$ Let $u:{\cal B}\mapsto {\cal D}$ be a  map for which  
 the diagram 
$${\cal B}\stackrel{u}{\longmapsto } {\cal D}$$
$$\theta \searrow    \swarrow \vartheta $$
$$X$$
commutes, so that $\theta (\xi ) = \vartheta (u(\xi ))$ for $\xi \in {\cal B}
.$ Let $Y$ be another Hausdorff space  and $\phi : Y \mapsto X$ 
be a continuous map. Let ${\cal B}^{\# }$ and ${\cal D}^{\# }$ be the 
retractions of ${\cal B}$ and ${\cal D}$ by $\phi $ respectively . 
Define the map 
$j^{\# } (u ):{\cal B}^{\# }{\mapsto } {\cal D}^{\# }$  by 
$$ j^{\# } (u )(y,\xi )=(y,u(\xi )).$$
Then 
$$\vartheta^{\# } ( j^{\# } (u ) ((y,\xi )))=\vartheta^{\# } (  (y,u(\xi ))=
 y =\theta^{\# }
(y,\xi ),$$
for $(y,\xi )\in {\cal B}^{\# }$, so that the diagram; 
$${\cal B}^{\# }\stackrel{j^{\# } (u )}{\longmapsto } {\cal D}^{\# }$$
$$\theta^{\# } \searrow  \quad \swarrow 
\vartheta^{\# } $$
$$Y$$
commutes.

Suppose $u:{\cal B}\mapsto {\cal D}$ is a continuous and open map.
It is clear that the map $j^{\# } (u )$ is the restriction of the map
$(j,u ): Y\times {\cal B}\mapsto Y\times {\cal D},$ where $j$ is the 
identity map from $Y$ to itself and $(j,u)(y,\xi ) = (y,u(\xi)).$ Clearly, 
$(j,u )$ is a continuous, open map. Let ${\tilde U}
\subset {\cal B}^{\# }$ be an open set. Then there exists an open set $U 
\subseteq Y\times {\cal B}$ such that ${\tilde U} = U\cap {\cal B}^{\# } .$
Let $j^{\# } (u ) {\tilde U} = {\tilde V}$ and $(j,u ) (U) =V.$ Now
$V$ is an open set in $Y\times {\cal D}$ and ${\tilde V} \subseteq V\cap 
{\cal D}^{\# }.$ Note that  
if $(y,\xi )\not \in {\cal B}^{\# },$ then $\phi (y) \not = \theta (\xi),$
and therefore $ \vartheta (u(\xi)) = \theta (\xi ) \not = \phi (y),$ 
which implies that $(y,u(\xi ))\not \in {\cal D}^{\# }.$ Therefore,
if  $x\in V\cap {\cal D}^{\# }
$ is the image of $z\in U,$ then $z$ cannot be outside of ${\cal B}^{\# }.$
 This implies that ${\tilde V} = V\cap {\cal D}^{\# },$ which shows that 
${\tilde V}$ is an open set in ${\cal D}^{\# }.$
 Hence $j^{\# } (u )$ is an open map. 

Now we turn to the construction of a particular type of  Banach \newline
\noindent (semi-)bundle.
 Let the Banach (semi-)bundle 
${\underline {\cal B}}=({\cal B}, \theta)$ over $X$ with ${\cal B}= {\cal H}
\times X$ be such that ${\cal H}$ is a Banach space, $X$ is a Hausdorff 
space and $\theta (\xi , x) = x.$  Suppose that there is an equivalence 
relation $R$ given on $X.$ Let $r$ be the canonical mapping from $X$ to 
$X/R.$ For $x\in X,$ let $r(x)\in X/R$ be the equivalence 
class to which $x$ belongs. Define ${\underline {\cal B}^R}=({\cal B}^R, 
\theta ^R)$ over $X/R$ by letting $ {\cal B}^R = {\cal H} \times {X/R}$ and 
$\theta^R (\xi ,r( x)) =r( x).$ Clearly, both bundles ${\underline {\cal B}}
$ and ${\underline {\cal B}^R}$ are trivial bundles with constant fibre 
${\cal H}.$ 

\begin{pro}
 The Banach 
bundle retraction $${\underline {\cal B}^R}^{\# }=({{\cal B}^R}^{\# }, 
{\theta ^R}^{\# })$$ of ${\underline {\cal B}^R}$ by $r$ is topologically 
equivalent to ${\underline {\cal B}}=({\cal B}, \theta).$
\end{pro}

\noindent Proof: The two Banach bundles ${\underline {\cal B}^R}^{\# }$
and ${\underline {\cal B}}$ have the same base space $X.$ 
\begin{eqnarray*} 
{{\cal B}^R}^{\# }&=&\{(x' ,(\xi , r(x)) ) : \theta (\xi , r(x))
= r(x' ), x',x\in X, \xi \in {\cal H} \}\\
&=&\{(x' ,(\xi , r(x)) ) : x' \in r(x),  x',x\in X, \xi \in {\cal H}\},
\end{eqnarray*}
and for $x\in X,$ ${\cal B}_x =\{(\xi , x): \xi \in {\cal H}\},$ while
${{\cal B}_x^R}^{\# } =\{(x,(\xi ,r(x))): \xi \in {\cal H}\}.$ Clearly, the
 mapping $(\xi , x)\mapsto (x,(\xi , r(x)))$ is a homeomorphism.

\begin{flushright}  \large$\diamond$\\
\end{flushright}

 A \sf cross-section \rm of ${\underline {\cal B}}$ is a function
 $f:X\mapsto
{\cal B}$ such that $f(x)\in {\cal B}_x $ for each $x\in X.$ 
The linear space of all  continuous cross-sections of ${\underline {\cal B}}$
is denoted by $ C({\underline {\cal B}})$ and the 
subspace of $ C({\underline {\cal B}})$ consisting of 
those cross-sections which vanish outside some compact set is denoted by
$ C_0({\underline {\cal B}}).$ The set of all bounded 
cross-sections is denoted by $B({\underline {\cal B}}).$

We say that ${\underline {\cal B}}$ has \sf enough continuous cross-sections
 \rm if for every $\xi \in {\cal B}$ there exists a continuous cross-section
 $f:X\mapsto {\cal B}$ for which $f(\theta (\xi ))=\xi.$

An unpublished result by A.Douady and L.dal Soglio-H\'erault about
 the existence of enough continuous cross-sections
states that if $X$ is either paracompact or locally compact, then 
 every Banach bundle over $X$ has enough continuous cross-sections
(see  Fell\cite{fel1}, p.324).

 Let $1\leq p<\infty .$ A cross-section of 
${\underline {\cal B}}$ is said to be \sf $p^{th}$-power summable \rm if it is locally
$\mu$-measurable and $$\| f \|_p = \bigl(\int_X \|f(x) \|^p d\mu (x) \bigr)^{1/p}
<\infty.$$
\noindent The space of all $p^{th}$-power summable cross-sections
is denoted by $L_p ({\underline {\cal B}};\mu ).$

 $L_p ({\underline {\cal B}};\mu )$ is a Banach space under the norm 
$\|\quad \|_p$ defined above.

The space $L_{\infty } ({\underline {\cal B}};\mu )$ is defined to
be the space of all $\mu$-essentially bounded cross-sections of 
${\underline {\cal B}}.$

 $L_{\infty } ({\underline {\cal B}};\mu )$ is a Banach space under the
norm $\|f\|_{\infty } = \mu {\hbox{-ess }}\sup_{x\in X} \|f(x)\|.$

  Let $Y$ be another locally compact Hausdorff space with a regular 
Borel measure $\nu.$
  Let $\kappa: X\times Y \mapsto X $ be the surjection $(x,y)\mapsto x$.
 Then the Banach \newline
\noindent (semi-)bundle  retraction ${\underline {\cal E}}=
({\cal E}, \rho)$ by $\kappa$ is a bundle over 
  $X\times Y$ whose bundle space ${\cal E}$ can be identified with ${\cal B}
 \times Y.$
  The bundle projection is given by $\rho :(\xi , y)\rightarrow  
(\theta (\xi ),y).$
 For each $x\in X,$  ${\underline {\cal E}}_{\{x\} \times Y}$ is the 
trivial bundle 
  with constant fibre ${\cal B}_x .$ Therefore, for a given $h\in 
 {\cal C}_0
  ({\underline {\cal E}})$ and for each $x$ in $X,$ the Bochner integral
 $\int_Y h(x,y) d\nu (y)$
   exists and  will belong to ${\cal B}_x .$

The following result has been  proved by Fell\cite{fel} 
for Banach bundles. The proof is similar in the context of Banach semi-bundles. 

\begin{lem}
 For each $h\in { {\cal C}_0}({\underline {\cal E}})$
the map $\ell (x)= \int_Y h(x,y) d\nu (y)$ is a continuous cross-section 
of the Banach semi-bundle  ${\underline {\cal B}}.$
\end{lem}

\subsection{ The p-induced representations of locally compact 
 groups and $L_p (\pi )$ spaces}

Let $G$ be a locally compact group and let $H$ be a closed subgroup of $G.$
Suppose that $\pi$ is a representation of $H$ on a Banach space
 ${\cal H}(\pi ).$
Let  $\mu$ be any quasi-invariant measure, in the homogeneous space $X=G/H$
of  right cosets, which belongs to a continuous $\rho$-function. 
For $1\leq p <\infty ,$ let us denote
 by $L_p (\pi ,\mu )$ the set of all functions  $f$ from G to 
a Banach Space $\cal {H}\rm (\pi)$ such that 

\par \noindent
(1) $\langle f(x),v\rangle $ is a Borel function of $x$ for all $v \in {\cal {H} \rm (\pi)}
^{*} ;$
\par \noindent (2) $f$ satisfies the  covariance condition $f(hx)= {\pi }_h 
f(x)$  for all $h\in H $ and $x\in G;$ and
\par \noindent
(3) $\| f \|_p =\biggl ({\int_{G\over H}  \|  f(x)
\|^p d\mu (z)}\biggr )^{1\over p} < \infty$.

Note that the integrand in the above integral is constant on each right coset $Hx$ and hence
defines a function on $X$. When functions equal almost everywhere are
identified, $L_p (\pi ,\mu )$ becomes a Banach space under the norm defined by
(3)(for which we use the same symbol $L_p (\pi ,\mu )$ ).

For each $x,y \in G $ and $f\in L_p (\pi ,\mu )$, let us define a mapping
${}^{\mu }U_y^ \pi $ on $L_p (\pi ,\mu )$ by
\begin{eqnarray}
\quad \quad \quad ({}^{\mu }U_y^ \pi f)(x)  := \lambda (x,y)^{1\over p} f(xy), 
\end{eqnarray} 
 where  $\lambda (\cdot ,y)$ is the Radon-Nikodym derivative of the measure 
$\mu_y$ with respect to the measure $\mu$. Then, it can be easily seen that
${}^{\mu }U^ \pi$ is a representation of the group $G$ on the Banach space
$L_p (\pi ,\mu ).$ Also, given two 
 quasi-invariant measures  $\mu $ and $\mu ^{'}$ 
 on $X,$ there exists an isometry $W$ from $L_p (\pi , \mu )$
onto $L_p (\pi , \mu^{'} )$ such that $W (^{\mu}U_y^{\pi})=
(^{\mu^{'}}U_y^{\pi}) W $ for all $y\in G.$ In other words, the two 
representations $ ^{\mu}U_y^{\pi}$ and $ ^{\mu^{'}}U_y^{\pi}$
are equivalent. (cf. Mackey\cite{mack1} ).

 The  equivalence
class of $^{\mu } U^{\pi}$ (denoted by $U^ \pi $) is called  the \sf 
representation
of $G$ induced by the representation $\pi $ of $H$.\rm The corresponding Banach 
space of (equivalence classes) of functions is denoted by $L_p (\pi )$. (The most
appropriate notation for the p-induced representation (induced by $\pi$) 
would be $U^{\pi }_p ;$ but for simplicity of notation we use $U^ \pi $ unless
the former is necessary to avoid confusion.) 

Let $\pi $ and $\gamma $ be representations of the locally compact group $G$.
A bounded linear operator $T$ from ${\cal H}(\pi )$ to ${\cal H}(\gamma)$ is called an 
\sf intertwining operator \rm for $\pi$ and $\gamma$ if $\pi (x) T 
= T\gamma (x) $ for all $x\in G.$
The vector space of all intertwining operators is denoted by $Int_G (\pi ,
\gamma )$ and the dimension (possibly infinite ) of this space, called the
\sf intertwining 
number\rm, is denoted by $\partial (\pi ,\gamma ).$

 Let $(\Omega ,\Sigma ,\mu )$  be a measure space. A Banach space $X$ is 
said to have the \sf Radon-Nikodym
property \rm with respect to $(\Omega ,\Sigma ,\mu )$ if for each $\mu$-continuous
vector  measure $F:\Sigma \rightarrow X$ of bounded variation there exists $g\in
L_1( X,\mu )$ such that $F(E) = \int_E g d\mu $ for all $E \in \Sigma .$
A Banach space $X$ has the Radon-Nikodym property if $X$ has the
Radon-Nikodym property with respect to every finite measure space
(see Gretsky and Uhl\cite{gre}, Chapter III).

 Let $(\Omega ,\Sigma ,\mu )$ be a $\sigma $-finite
measure space, $1\leq p <\infty ,$ and let $X$ be a Banach space. It is well known that
$L_p(\Omega , X,\mu )^* = L_{p'} (\Omega , X^*,\mu  ),$ where $1/p +1/{p'} =1,$ if and only if $X^*$ has the
Radon-Nikodym property with respect to $\mu.$ 
 Also,  if $(\Omega ,\Sigma ,\mu
)$ is a nonatomic finite measure space, then it can be seen that
 $L_p(\Omega , X,\mu )$ has the Radon-Nikodym property if and only if
 $1<p <\infty $ and
$X$ has the Radon-Nikodym property.

\sf Throughout our work, we assume that the Banach space ${\cal H }(\pi ) $ of
a representation $\pi $ of a subgroup $H$ of a group $G$ stays within the class
of spaces satisfying the Radon-Nikodym property. \rm 

 Let $\pi $ be a representation of a group $G$ on a Banach space
 ${\cal H}(\pi ).$  We define the map
 $\pi^* :G \mapsto U 
\bigl(({\cal H}(\pi ))^*\bigr)$ by letting
 $\pi^* (x) = (\pi (x^{-1}))^*.$ It can be easily seen that that $\pi^* $ 
is a representation 
of $G$ on the Banach space ${\cal H}(\pi^* )= ({\cal H}(\pi ))^* ,$
when ${\cal H} (\pi )$ is reflexive.
 Assume now that the Banach space 
${\cal H }(\pi ) $ is reflexive.
 Let us consider the 
Banach space $L_{p'} (\pi^*)$ and the induced representation $U_{p'}^{\pi^*}$ of $G.$
 The dual pairing between $L_p (\pi)$
 and  $L_{p'} (\pi^*)$  is given by  $$\langle f,g\rangle =\int_{G\over H} 
\langle f(x),g(x)\rangle 
 d\mu (x), 
{\hbox { for }}  f \in L_p (\pi){\hbox { and }}
 g \in L_{{p}'} (\pi^*).$$ 
The above integral is well defined since, 
for any $h\in H$ and $x\in G,$
 \begin{eqnarray*}
\langle f(hx),g(hx)\rangle &=&\langle \pi (h)f(x),\pi^* (h)g(x)\rangle ,\\
&=& \langle \pi (h)f(x),(\pi (h^{-1}))^* g(x)\rangle ,\\
&=&\langle f(x),g(x)\rangle .\\
\end{eqnarray*}
\noindent Also, for any $y\in G$,
\begin{eqnarray*}
\langle U_p^\pi (y) f,U_{p'}^{\pi^*} (y)g\rangle  &
=&\int_{G\over H}\langle \lambda (x,y)^{1\over p} f(xy),
\lambda (x,y)^{1\over {p'}} g(xy)\rangle d\mu (x)
\\
&=& \int_{G\over H}\lambda (x,y) \langle f(xy), g(xy)\rangle d\mu (x)\\
&=&\langle f,g\rangle ,
\end{eqnarray*}
the last equality of which was obtained
 by changing variables $x\mapsto xy.$ This implies that
\begin{eqnarray}
U_{p'}^{\pi ^*}(y) =(U_p^\pi  (y^{-1}))^* & = & (U_{p}^{\pi })^* (y),
 {\hbox { for all } } y\in G.
\end{eqnarray}

Let $1\leq p <\infty .$ Let us define a convolution $g*f $ for $g\in L_p (\pi )$ 
and $ f\in L_1 (G),$ 
by 
$$(g*f)(x) := \int_G (\lambda_H (x,y^{-1}))^{1\over{p}}g(xy^{-1})f(y)d\nu_G (y).$$
It is not difficult to prove that $g*f $ belongs to $L_p (\pi ),$ 
$g*(h*f)=  (g*h)*f $
for all  $g\in L_p (\pi )$ and $h,f\in L_1 (G)$
and that $ L_p (\pi )$ is an $ L_1 (G)$-module.

Let  $Hom_G (L_p (\pi ), {L}_{q} (\gamma ))$ denote the Banach space of all continuous $G$-module 
homomorphisms from $L_p (\pi )$ to $ {L}_{q} (\gamma )$ (Rieffel\cite{rief2}). 

\begin{pro}
\begin{eqnarray}
Hom_G (L_p (\pi ), {L}_{q} (\gamma )) &= & Int_G (U_p^{\pi }, U_q^{\gamma }).
\end{eqnarray}
\end{pro}

\noindent
Proof:
Let $T$ be any bounded linear operator from $L_p (\pi )$ to $L_{q}
 (\gamma )$ and $T^* $ be its adjoint operator. For any $g\in L_p(\pi ), f\in L_1 (G)$ and $k\in L_{q'}
 (\gamma^* ),$ 
\begin{eqnarray*}
\langle T(g*f) , k\rangle & = & \langle g*f , T^* k\rangle ,\\
 & = & \int_{G\over H}\langle \int_G (\lambda (x,y^{-1}))^{1\over{p}} g(xy^{-1})f(y)d\nu_G (y)
 ,T^* k(x)\rangle d\mu_H (x)  ,\\
  & = & \int_G f(y)\int_{G\over H}\langle (U_p^{\pi }{(y^{-1})}g)(x) , 
T^* k(x)\rangle d\mu_H (x)
 d\nu_G (y) ,\\      
 & = & \int_G f(y)\langle U_p^{\pi }{(y^{-1})}g , T^* k\rangle d\nu_G (y).
\end{eqnarray*}
\noindent
Hence, 
\begin{eqnarray}
\langle T(g*f) , k\rangle & = & \int_G f(y)\langle TU_p^{\pi }
{(y^{-1})}g , k\rangle d\nu_G (y). 
\end{eqnarray}
\noindent
On the other hand,
\begin{eqnarray*}
\langle T(g)*f , k\rangle  & = & \int_{G\over K} \langle 
(T(g)*f)(x) , k(x)\rangle d\mu_K (x ),\\
& = & \int_{G\over K} \langle \int_G (\lambda (x,y^{-1}))^{1\over{q'}}T(g)(xy^{-1})f(y) ,
 k(x)\rangle d\mu_K (x ) d\nu_G (y),\\
& = & \int_G f(y) \int_{G\over K} \langle (U_q^{\gamma }
{(y^{-1})}Tg)(x) , k(x)\rangle d\mu_K (x) 
d\nu_G (y). 
\end{eqnarray*}
\noindent
Therefore, 
\begin{eqnarray}
 \langle T(g)*f , k\rangle & = & \int_G f(y) \langle U_q^{\gamma }
{(y^{-1})}Tg , k\rangle  d\nu_G (y).
 \end{eqnarray}

\noindent
If $T\in Hom_G (L_p (\pi ), L_{q}(\gamma )),$ we see, by (5) and (6), that
\begin{eqnarray}
TU_p^{\pi }(y)g &=& U_q^{\gamma }(y)Tg ,
\end{eqnarray}
for almost all $y\in G.$ 
By continuity, (7) is true for all $y\in G.$  
Hence \newline \noindent $T\in Int_G (U_p^{\pi }, U_q^{\gamma })$. Conversely, 
 $T\in Int_G (U_p^{\pi }, U_q^{\gamma })$ implies \newline \noindent
 $T\in Hom_G (L_p 
(\pi ), L_{q}(\gamma )),$ by (5) and (6). Hence, (4) follows.

\begin{flushright}  \large$\diamond$\\
\end{flushright}

\section{Some important results on $\lambda$-functions}

First, we intend to prove an integral formula
which involves integration on coset spaces. Secondly,  
 the notion of \sf disintegration  of measures \rm (which 
has been discussed in a
number of places in the literature (see, for example, Mackey\cite{mack1}, 
Halmos\cite{hal1})) is dealt with. Here, we derive an identity among 
$\lambda $-functions of a particular set of 
subgroups of a given group. 

\begin{lem} Let $G$ be a locally compact group. Let $H$ and $ K$ be
 subgroups of
 $G$ with $K\subseteq H.$ Then there exist
 positive quasi-invariant measures $\mu_{ K}$ on $G/{ K}$, 
$\mu_{H}$ on $G/{H}$ and 
${\tilde{\mu }}$ 
on $H/{ K}$   such that, for $F\in {C_0}(G/K),$

\begin{eqnarray}
\int_{G\over { K}} F(z) d\mu_{K}(z) & = & \int_{G\over {H}}
\Biggl (\int_{H\over{K}}
 {\lambda_{ K}(y,t)\over {\lambda_{H}(y,t)}}F(yt) d{\tilde{\mu
}}(y)\Biggr ) d\mu_{H}(t),
\end{eqnarray}
 whenever the integrals exist.
\end{lem}

\noindent Proof: As discussed in Corollary 1.2
(see also Reiter\cite{rei}, p.158, Mackey\cite{mack1}), there exists a 
continuous, strictly positive function 
${\rho_{ K}}$ on $G$ and a positive measure $\mu_{ K}$ on 
$G/{ K}$ such that
\begin{eqnarray}
\int_G f(u) d\nu_{G}(u) & = & \int_{G\over{ K}}\Bigl(\int_{K}{1\over 
{\rho_{ K}(sz)}} f(sz) d\nu_{ K} (s) \Bigl )d\mu_{K}(z),
\end{eqnarray}
  for $ f\in {C_0}(G).$ 
  
Also, by the same reasoning, there  exists a continuous, strictly positive function 
$\rho_{H}$ on $G$ and a positive measure $\mu_{H}$ on $G/{H}$ such that
$$\int_G f(u)d \nu_G (u) = \int_{G\over {H}}\Bigl(\int_{H}{1\over {\rho_H (ht)}} 
f(ht) d\nu_H (h) \Bigl )d\mu_{H}(t) .$$ 

Let $\tilde{\rho } = {\rho_{ K} / \rho_{H}}.$ We see that 
$$\tilde{\rho }(sx) = {\rho_{K}(sx)
 / {\rho_{H}(sx)}} = ({\Delta_{ K} (s)}/ \Delta_{H} (s))
\tilde{\rho }(x),$$ 
for $s\in { K}$ and $x\in G. $ Thus
$\tilde{\rho },$ restricted to $H,$ is a $\rho$-function for the
homogeneous space ${H/ { K}}.$ If we let $\tilde\mu $ be a 
quasi-invariant measure associated with this $\rho $-function, we have
\begin{eqnarray}
\int_G f(u) d\nu_G (u) & = & \int_{G\over{H}}\int_{H\over{K}} 
\Bigl(\int_{ K}
{\rho_H (sy) \over {\rho_{K}(sy)\rho_{H}(syt)}} f(syt) d\nu_{K} (s) \Bigl )
d{\tilde{\mu }}(y) \mu_{H}(t).    
\end{eqnarray}

 By Reiter\cite{rei}, p.165, for a given $F\in {C}(G/K),$ 
there exists a function \newline \noindent $f\in {C}(G)$ such that
\begin{eqnarray}
F(\dot{z}) & = & \int_{K}{1\over {\rho_{K}(sz)}} f(sz) 
d\nu_{K} (s),
\end{eqnarray} 

\noindent where $\dot{z} = p_K (z).$
Comparing equations (9) and (10), and using (11),
we see that
\begin{eqnarray*}
\int_{G\over {K}} F(z) d\mu_{K}(z) & = & \int_{G\over {H}}\int_{H\over {K}}
{{\rho_{K}(yt) \rho_H (y) }\over {\rho_{K}(y)\rho_{H}(yt)}}F(yt)
 d\tilde \mu (y) \mu_{H}(t)\\
 &= &\int_{G\over {H}} \Biggl(\int_{H\over {K}}{\lambda_{K}(y,t)
\over {\lambda_{H}(y,t)}}F(yt) d{\tilde {\mu }}(y)\Biggr ) d\mu_{H}(t), 
 \end{eqnarray*}
for any $F\in {C_0}(G/K),$ and (8) is proved.

\begin{flushright}  \large$\diamond$\\
\end{flushright}

Let $\mu_{H}$ be a given quasi-invariant measure on $G/H$ with the corresponding 
$\lambda$-function $\lambda_{H}.$ Consider the homeomorphism $\phi_x : G/H^x \mapsto G/H $ given  by
$\phi_x (u) = xu .$ Define  
a measure $\mu_{H^x}$  on $G/H^x$ by  $\mu_{H^x}(E) = \mu_H  ( \phi_x (E))$
whenever $E$ is such that $x.E$ is measurable.
 Clearly, $\mu_{H^x}$ is quasi-invariant if and only if $\mu_H $
is. The corresponding $\lambda$-function of $\mu_{H^x}$ is denoted by 
$\lambda_{H^x}.$ Then, it can be easily seen that,
for $x,t\in G$ and for almost all $v\in G/H,$ 
\begin{eqnarray}
\lambda_{H^x} (x^{-1}v,t) & = & \lambda_{H} (v,t).
\end{eqnarray} 
which states the relationship between 
$\lambda_{H}$ and $\lambda_{H^x }.$

Let $\Delta = \{ (x,x):x\in G \}$ be the diagonal subgroup of $G\times G.$
Consider the right action of $\Delta$ on the coset space ${({G\times G})/ 
({H\times K})}.$ The stabilizer of the coset $(Hx,Ky)$ is $(H\times K)^{(x,y)}
\cap \Delta$ and the orbit is the double coset $(H\times K)(x,y)\Delta .$
Let  $\Upsilon $ be the set of all  double cosets $(H\times K) : \Delta $ 
of $G\times G;$ that is, the set of all orbits. For each 
$(x,y) \in G\times G,$ let $k(x,y)$ denote 
the $(H\times K) : \Delta $ double coset to which $(x,y) $ belongs. If 
$\nu_0 $ is any finite measure in $G\times G$ with the same null sets as
 Haar measure we define a measure  $\mu_{(H,K)} $  on $\Upsilon $ by 
$\mu_{(H,K)}(F) = \nu_0 (k^{-1} (F)) $ whenever $F$ is such that 
$k^{-1} (F)$ is measurable. Using Mackey's terminology, we call such a
 measure an \sf admissible measure \rm in $\Upsilon .$ We obtain the
 following result as a consequence of Lemma 11.1, Mackey\cite{mack1}. 

\begin{lem}
 Suppose that $H$ and $K$ are
regularly related(see  Mackey\cite{mack1}). Let $\Delta $ be the diagonal subgroup of $G\times G$ and $\Upsilon $
 denote the set of all $(H\times K):{\Delta }$ double cosets in $G \times G$. 
Then for each double coset $D(x,y)= H\times K (x,y) \Delta $ there exists a 
 quasi-invariant measure $\mu_{x,y}$ on $G /
 (H^x\cap K^y )$, $x,y\in G$, and $
\lambda_{H^x \cap K^y} $
with $\mu_{x,y}\succ \lambda_{H^x \cap K^y}$ such that
\begin{eqnarray}
\lambda_H (xts^{-1},s)\lambda_K (yts^{-1},s)\lambda_{H^x \cap K^y} (\underline{t} ,
\underline{s}^{-1}) &=& 1 ,
\end{eqnarray}
 for all $s,t\in G,$   and for almost all $(x,y)\in 
{(G\times G)}/{(H \times K)}.$ Moreover, $\lambda_{H^x \cap K^y} (t,s)$ is 
defined everywhere and continuous on $(G/(H^x \cap K^y )) \times G .$

 \end{lem}

\par \noindent 
Proof: Choose  two quasi-invariant measures $\mu_H $ and $\mu_K$ on $G/H$ and
$G/K$ respectively, which correspond to two continuous $\rho$-functions. 
Define a measure $\mu_{H\times K}$ in ${(G\times G)/(H\times K)}$ by 
$\mu_{H\times K}= \mu_H \times \mu_K $
 (see, for example, Halmos\cite{hal2}, p.144). Obviously, 
$\mu_{H\times K}$ is quasi-invariant to the action of $\Delta.$
  Let $\nu_0 $ be the 
 measure in ${(G\times G)}$ defined
 by $ \nu_0 (p_{H\times K}^{-1} (F))=\mu_{H\times K} (F).$ 
Let $\mu_{H,K} $ be an admissible measure in $\Upsilon $ corresponding to 
 $\nu_0 .$

Let $f$ be a function defined on $(G/H )\times (G/K) $. 
Suppose \newline \noindent $\int_{G\over H} \int_{G\over K} f(x,y) d\mu_H (x) d\mu_K (y)$ is 
 integrable. Changing the variables $x\mapsto xs$ and 
$y\mapsto ys,$ we get
\begin{eqnarray*}
\lefteqn{\int_{G\over{H}} \int_{G\over{K}} f(x,y) d\mu_H (x) d\mu_K (y)   }
      \\
& = & \int_{G\over{H}} \int_{G\over{K}}\lambda_H (x,s)\lambda_K (y,s)f(xs,ys) 
d\mu_H (x) d\mu_K (y)  \\
& = & \int_{{G\times G}\over{H \times K}} \lambda_H (x,s) \lambda_K (y,s)
f(xs,ys)d\mu_{H\times K} (x,y).
\end{eqnarray*}

\noindent
For each $(x,y) $ in ${(G\times G)}/{(H \times K)}$ let $r(x,y) = 
k(p_{H\times K}^{-1} (x,y)).$
If $H$ and $K$ are regularly related then $r$ defines a measurable equivalence 
relation(see Mackey\cite{mack1}). Then, by Lemma 11.1, Mackey\cite{mack1}, 
$\mu_{H\times K}$ is an 
integral of measures $\mu_{x,y} $, where $D(x,y)\in 
\Upsilon $, with respect to the measure $\mu_{H,K} $ in $\Upsilon $.
By Lemma 11.5, Mackey\cite{mack1}, each $\mu_{x,y} $ is a quasi-invariant measure  
on the orbit $r^{-1}(D(x,y)).$ 
 Using this disintegration, we have 
\begin{eqnarray*}
\lefteqn{\int_{{G\times G} \over {H \times K}} \lambda_H (x,s) 
\lambda_K (y,s)
f(xs,ys)d\mu_{H\times K}(x,y)    }
        \\
& =&\int_{D\in \Upsilon } \int_{\underline{t} \in {\Delta \over 
{{(H\times K)}^{(x,y)}\cap \Delta}}}\lambda_H (xt,s)\lambda_K (yt,s)
f(xts,yts) d\mu_{x,y} (\underline{t} ) d\mu_{H,K} (D), 
\end{eqnarray*}
\noindent
where $(x,y)$ is the coset representative of the coset $D(x,y)$. Identifying
the space ${\Delta / 
({{(H\times K)}^{(x,y)}\cap \Delta})}$ with $G/ (H^x \cap K^y)$ we can regard
$\mu_{x,y} $ as a measure on  $G/ (H^x \cap K^y)$. Then we have
\begin{eqnarray*}
\lefteqn{\int_{{G\times G} \over {H \times K}} \lambda_H (x,s) 
\lambda_K (y,s)
f(xs,ys)d\mu_{H\times K}(x,y)    }
        \\
& =&\int_{D\in \Upsilon } \int_{t \in {G\over {H^x \cap K^y}}}\lambda_H (xt,s)\lambda_K (yt,s)
f(xts,yts) d\mu_{x,y} ({t} ) d\mu_{H,K} (D), 
\end{eqnarray*}
 Changing variables ${t} \mapsto {ts}^{-1}$, in the 
integral on the right-hand side, we get
\begin{eqnarray}
\lefteqn{\int_{{G\times G} \over {H \times K}} \lambda_H (x,s) \lambda_K (y,s)
f(xs,ys)d\mu_{H\times K}(x,y)     }
\nonumber          \\
&= &\int_{D\in \Upsilon } \int_{t \in {G\over {H^x \cap K^y}}}
\lambda_H (xts^{-1},s) \lambda_K (yts^{-1},s)f(xt,yt) \nonumber \\
& &\quad \quad \quad \quad \quad \quad \quad \quad \quad \quad 
\quad   
 \lambda_{H^x \cap K^y} ({t} ,{s}^{-1})
 d\mu_{x,y} ({t} ) d\mu_{H,K} (D). 
\end{eqnarray}

\noindent 
On the other hand, if we start with ${\int \int }_{{G\times G}\over 
{H \times K}}f(x,y) d\mu_{H\times K}(x,y)$
and use Lemma Lemma 11.1, Mackey\cite{mack1}, we have
\begin{eqnarray}
\lefteqn{{\int \int }_{{G\times G}\over 
{H \times K}}f(x,y) d\mu_{H\times K}(x,y)  } \hspace{0.6in}    \nonumber  \\
&= & \int_{D \in \Upsilon } 
\int_{\underline{t} \in {\Delta \over{{(H\times K)^{(x,y)} \cap \Delta}}} }
f(xt,yt) d\mu_{x,y}(\underline{t})
 d\mu_{(H,K)} (D), \nonumber \\
&= & \int_{D \in \Upsilon } 
\int_{t \in {G \over{H^x \cap K^y} }}
f(xt,yt) d\mu_{x,y}({t})
 d\mu_{(H,K)} (D). 
\end{eqnarray}
\noindent
Hence from (14) and (15) we have  
\begin{eqnarray*}
\lambda_H (xts^{-1},s)
 \lambda_K (yts^{-1},s)\lambda_{H^x \cap K^y} ({t} ,{s}^{-1}) & = &1,
\end{eqnarray*}
 for all $ s \in G,$ for almost all $t\in {G /({H^x \cap K^y}) }$ and 
for almost all \newline \noindent $(x,y)\in 
{(G\times G)}/{(H \times K)}.$ For each such 
$(x_0,y_0)\in 
{(G\times G)}/{(H \times K)},$
\begin{eqnarray}
\lambda_H (x_0ts^{-1},s)
 \lambda_K (y_0ts^{-1},s)\lambda_{H^{x_0} \cap K^{y_0}} 
({t} ,{s}^{-1}) & = &1.
\end{eqnarray}
By continuity of $\lambda_H $ and $\lambda_K ,$ we see that (16) is true for all
$t\in {G /({H^{x_0} \cap K^{y_0}}) }.$ Furthermore, (16) implies that 
$\lambda_{H^x \cap K^y} ({t} ,{s} )$ is defined everywhere and continuous 
on $({G /({H^x \cap K^y}) })\times G,$
 which proves the Lemma.
 
\par \begin{flushright}  \large$\diamond$\\
\end{flushright}

The following result is a consequence of Lemma 2.2.
\begin{cor} Let $(x,y)\in G\times G$ such that the identity (13) holds.
Then for $ s\in H^x\cap K^y $ ,
\begin{eqnarray}
{{\Delta_H (h) \Delta_K (k)}\over {\Delta_G (s) \Delta_{H^x\cap K^y} (s)} } & =&
1,\quad \quad \quad \quad 
\end{eqnarray}
where $h=xsx^{-1}$ and $k=ysy^{-1} .$
\end{cor}

\noindent Proof: Let $t=s$ in the identity (13). Then we have
\begin{eqnarray}
\lambda_H (x,s)
 \lambda_K (y,s)\lambda_{H^x \cap K^y} ({s} ,{s}^{-1})& =&1 .
\end{eqnarray}
By property (a) of $\lambda$-functions given in Section 1.1, page 3, 
 this simplifies to
\begin{eqnarray}
\lambda_H (x,s)
 \lambda_K (y,s)& =&\lambda_{H^x \cap K^y} (e ,{s}). 
\end{eqnarray}
Consider $s\in H^x\cap  K^y $. Then  $s= x^{-1} h x =y^{-1} ky$  for some 
$h\in H$ and $k\in K.$ For such an $s,$ we have by properties (a) and (b) 
of $\lambda$-functions, page 3,
\begin{eqnarray}
\lambda_H (x,s) = \lambda_H (x,x^{-1} h x)  & =& {{\Delta_H (h)} \over {\Delta_G (h)}},
\end{eqnarray}
 Similarly, 
\begin{eqnarray}
\lambda_K (y,s)& =&{{\Delta_K (k)} \over {\Delta_G (k)}},{\hbox{ and }}
\lambda_{H^x \cap K^y} (e ,{s}) ={{\Delta_{H^x \cap K^y} (s)}
 \over {\Delta_G (s)}}.
\end{eqnarray}
\noindent Using (19),(20) and (21), we obtain
\begin{eqnarray}
{{\Delta_H (h)} \over {\Delta_G (h)}} {{\Delta_K (k)} \over {\Delta_G (k)}}
 &=&
{{\Delta_{H^x \cap K^y} (s)}
 \over {\Delta_G (s)}}.
\end{eqnarray}
But $\Delta_G (h) = \Delta_G (x^{-1} h x)=\Delta_G (s) =\Delta_G ( y^{-1} k y)
=\Delta_G (k),$ hence (22) simplifies to
\begin{eqnarray}
{{\Delta_H (h) \Delta_K (k)}\over {\Delta_G (s) \Delta_{H^x\cap K^y}}(s) } & =&
1,\quad \quad \quad \quad 
\end{eqnarray}
as required.

\section{Projective tensor products and $A_p^q$ spaces}

\subsection{Construction of the convolution formula}

Let $G$ be a second countable locally 
compact group, with  closed subgroups $H$ and $K.$ Thus, the corresponding 
homogeneous spaces are Hausdorff and second countable, which in turn implies 
that any Borel measure on such spaces is regular. In addition, we will assume 
that $H$ and $K$ are regularly related (\cite{mack2}). $\mu_H $ and $\mu_K$ 
will denote fixed quasi-invariant measures on $G/H$ and $G/K,$ respectively.
We choose a family of quasi-invariant measures $\{ \mu_{x,y} : x\in G/H ,
 y\in G/K \},$ where $\mu_{x,y} $ is a measure on $G/(H^x \cap K^y ) ,$
in such a manner that for a function $f$ defined and integrable on $(G/H)\times 
(G/K )$, we have 
$$\int_{G/H} \int_{G/K} f(x,y) d\mu_H (x) d\mu_K (y) = \int_{D(x,y)\in \Upsilon }
\int_{t\in {G\over {H^x \cap K^y}}} f(xt,yt) d\mu_{x,y} (t) d\mu_{H,K} (D),$$
by disintegration of measures discussed in Lemma 2.2. For a given 
$\mu_{x,y},$ $\rho_{H^x \cap K^y}$ and $\lambda_{H^x \cap K^y}$ will denote
 the corresponding $\rho$-function and the 
$\lambda $-function respectively. For any $x\in G ,$
the quasi-invariant measure $\mu_{H^x} $ on $G/H^x $ will always considered to be
$\mu_{H^x} = \mu_{H} \circ \phi_x ,$ where $\phi_x :G/H^x \mapsto G/H $ 
is the homeomorphism given by $\phi_x (u) =xu .$ By $
\rho_{H^x} $ we mean the corresponding $\rho$-function of the above
 $\mu_{H^x}.$

 $\pi $ and $\gamma $ will denote representations of $H$ and $K$ 
on Banach spaces ${\cal H}(\pi )$ and 
${\cal H }(\gamma ),$ respectively.

Let $L_p (\pi) \otimes^{\sigma }  L_{q'} (\gamma^*)$ denote the projective tensor product
(\cite{gro}) of $L_p (\pi)$ and $ L_{q'} (\gamma^*)$ as Banach spaces so that ${\sigma }$
is the greatest cross-norm.
 Let $L$ be the closed linear subspace of
  $L_p (\pi) \otimes^{\sigma }  L_{q'} (\gamma^*)$ which is spanned by all 
the elements of the form
 $$  U_p^{\pi } (s) f\otimes g - f\otimes (U_{q}^{\gamma })^* (s)
g  ,\quad  s\in G, f\in L_p (\pi), g\in L_{q'}.$$ 
The quotient Banach space $( L_p (\pi) \otimes^{\sigma }  L_{q'} (\gamma^*) )/L$ is called 
the \sf $G$-module tensor product\rm, and is denoted by  $L_p (\pi) \otimes_G^{\sigma } 
 L_{q'} (\gamma^*) .$  
 Then we have a natural isometric isomorphism 
\begin{eqnarray} 
Int_G (U_p^{\pi } , U_q^{\gamma }) &\cong & (L_p (\pi) \otimes_G^{\sigma }  L_{q'} (\gamma^*)) ^* 
\end{eqnarray}
(see \cite{rief2}, 2.12 and 2.13), and the ultraweak*-topology on $Int_G (U_p^{\pi } , U_q^{\gamma })$
corresponds to the weak*-topology on $(L_p (\pi) \otimes_G^{\sigma }  L_{q'} (\gamma^*)) ^* $ 
(\cite{rief1}, Theorem 1.4).

Recall  that the space $A_p^q$ in the classical case consists of convolutions
of complex-valued functions of  $L_p (G)$ and $L_q (G)$ (see, for example, Rieffel\cite{rie1}).
 Our aim is to construct
$A_p^q$ spaces using spaces of induced representations, $L_p (\pi )$ and $L
_{q' } (\gamma^* ),$ which are spaces of vector-valued functions. Therefore, our 
 task is to construct a formula (Definition 4.7) for a convolution  of functions 
in 
$L_p (\pi )$ and $L_{q' } (\gamma^* ).$ The case where $G/H$ and $G/K$
are not compact is similar to that in the classical case (see H\"ormander\cite{Ho})
 in the sense that 
the non-triviality of the  tensor product
$ L_p (\pi )\otimes_G^{\sigma} L_{q'} (\gamma^* ) $ depends on the value of 
$1/p +1/q' $ as the following theorem states.   

\begin{thm}  Let ${1/ p }+{1/{q' }}<1,$ $1< p,q'<\infty.$
 Suppose that for any given compact set $F$ in $G,$ there 
exists $x\in G$ such that $HFx\cap HF =\emptyset$ and  
$KFx\cap KF =\emptyset.$
    Then $$ L_p (\pi )\otimes_G^{\sigma} L_{q'} (\gamma^* ) = \{{\underline 0}\}.$$
\end{thm}

We do not know whether Theorem 5.1 is true in the absence of 
the condition that there exists an element $x\in G$ such that $HFx\cap HF 
=\emptyset$ and  $KFx\cap KF =\emptyset$ for a given compact set $F$ in $G.$

Let us turn to the construction of the convolution formula. The following 
proposition states a result 
that equips us with the necessary ground work.

\par
\noindent
\begin{pro} Let $1\leq p,{q'} <\infty .$
For $\sum_{i=1}^{\infty } f_i \otimes g_i $ in $L_p(\pi ) \otimes^{\sigma} L_{q'} 
(\gamma^* )$ and for almost all $x\in {G/ H}$ and $y\in {G/ K},$ 
       $$ \sum_{i=1}^{\infty } f_i (x) \otimes g_i (y) \in {\cal H}(\pi )
\otimes^{\sigma }{\cal H}(\gamma^* ).$$
 \end{pro}

Our objective is to define a mapping on $L_p(\pi ) \otimes^{\sigma} L_{q'} 
(\gamma^* )$ so that its image space is a generalisation of the space of convolutions
as in the classical case.
Let us consider the integral
\begin{eqnarray}
\quad \quad  \int_G f(xt)\otimes^{\sigma }g(yt) d\nu_G (t) ,
\end{eqnarray}
 where $f\in L_p (\pi ) ,g\in L_{q'} (\gamma^* )$ and 
$(x,y)\in G\times G.$ It is easy to see that the norm of the integrand is
 constant on the subgroup 
$H^x \cap K^y $ of $G;$ for, if $t=x^{-1} h x =y^{-1} k y$ for some $h\in H$
and $k\in K,$ then $f(xt)\otimes g(yt) = \pi (h) f(x)\otimes \gamma^* (k) g(y)
.$ This implies that the space over which we integrate must reduce to 
${G/( {H^x\cap K^y})},$ in order to avoid the integrand becoming too large.
The integrand is constant over a given coset of ${G/( {H^x\cap K^y})}$ 
if 
 $$f(xst) \otimes g(yst) = f(xt) \otimes g(yt),$$
for all $s\in {H^x\cap K^y} .$ But 
$$f(xst) \otimes g(yst) ={\pi}^x (s) f(xt) \otimes  
{\gamma^*}^y (s) g(yt) .$$ 
This suggests that the integrand must have its value at $(x,y)$ 
in the quotient space  of ${\cal H}(\pi )
\otimes^{\sigma }{\cal H}(\gamma^*)  ,$ in which we have the equality 
$${\pi}^x (s) f(xt) \otimes  {\gamma^*}^y (s) g(yt)= f(xt) \otimes g(yt).$$
This calls for the following
definition. 

\begin{dfn} 
For any $x,y\in G,$ the \sf subspace $ {\cal H}_{x,y}$ \it of ${\cal H}(\pi )
\otimes^{\sigma }{\cal H}(\gamma^*) $ is defined to be the closed linear
 span of  elements of the form  
$${{\pi }}^x (b) 
 \xi \otimes \eta - \xi \otimes ({{\gamma}}^y (b))^* \eta,$$
 where $b \in H^x \cap K^y , \xi \in {\cal H}(\pi )$ and $ \eta \in {\cal H}
(\gamma^*). $
The quotient Banach space \newline 
\noindent ${\cal H}(\pi )\otimes
 ^{\sigma }{\cal H}(\gamma^*) / {\cal H}_{x,y }$ is denoted by 
${\cal A}_{x,y}.$ 
\end{dfn}
Note that, using the notation in Rieffel\cite{rief2}, ${\cal A}_{x,y}$
can be written as \newline \noindent ${\cal H}(\pi^x )\otimes_{H^x \cap K^y}^{\sigma }
{\cal H}({\gamma^*}^y).$
\begin{pro} 
 The spaces $\{{\cal H}_{x,y} : x,y\in G\},$ and hence the spaces \newline
\noindent $\{{\cal A}_{x,y} : x,y\in G\},$ satisfy the property  that 
\begin{eqnarray}
\quad \quad \quad \quad \quad \quad \quad {\cal H}_{xs,ys}= {\cal H}_{x,y} &{\hbox { and }} &{\cal A}_{xs,ys}=
 {\cal A}_{x,y},
\end{eqnarray}
 for any $s\in G.$
 \end{pro}

\noindent 
Proof: For $s\in G,$ the space ${\cal H}_{xs,ys}$ is the closed linear span 
of  elements of the form $${\pi }^{xs} (b) 
 \xi \otimes \eta - \xi \otimes ({\gamma}^{y s}(b))^* \eta ,$$ where 
$b \in H^{xs} \cap K^{ys} , \xi \in H(\pi )$ and $ \eta \in H(\gamma )$
 with 
$${{\pi }}^{xs} (b) = {\pi }(xsbs^{-1}x^{-1})= {{\pi }}^{x} (sbs^{-1}).$$
Since $b \in H^{xs} \cap K^{ys},$ there exist $h\in H $ and $k\in K $
such that \newline
\noindent $ b = s^{-1} x^{-1} h xs = s^{-1} y^{-1}
 kys .$ Hence  $sbs^{-1} = x^{-1} h x= y^{-1} ky,$ showing
that \newline
\noindent $sbs^{-1}\in H^x\cap K^y .$ 
Therefore, 
$${\pi }^{xs} (b) \xi \otimes \eta - \xi \otimes ({\gamma}^{y s}(b))^* 
\eta ={\pi }^{x} (sbs^{-1}) 
 \xi \otimes \eta - \xi \otimes ({\gamma}^{y }(sbs^{-1}))^* \eta ,$$
with $sbs^{-1} \in H^{x} \cap K^{y} , \xi \in H(\pi )$ and $
 \eta \in H(\gamma ).$
This implies that ${\cal H}_{xs,ys}\subseteq {\cal H}_{x,y},$
 \noindent
 which in turn
gives us that ${\cal H}_{x,y}={\cal H}_{xss^{-1},yss^{-1}} \subseteq 
{\cal H}_{xs,ys},$ for all $s\in G.$ Hence (28) follows.

\begin{flushright} \large$\diamond$\\
\end{flushright}

For $u\otimes v\in {\cal H}(\pi )\otimes^{\sigma }{\cal H}(\gamma^*),$ we use  
the notation $u\otimes_{x,y} v$ to denote the element of ${\cal A}_{x,y}$ 
to which $u\otimes v$ belongs. Then the integral (27) must be written in the
 form
\begin{eqnarray}
\int_{G\over {H^x\cap K^y}}
 f(xt)\otimes_{x,y }g(yt) d\mu_{x,y} (t), 
\end{eqnarray}
\noindent for a suitably chosen quasi-invariant measure $\mu_{x,y}$ on 
the homogeneous space ${G/( {H^x\cap K^y})}.$ For each $x,y\in G,$
the value of the 
integral belongs to the quotient Banach space ${\cal A}_{x,y}.$ The next obvious
 step in this construction is to check whether the integral is finite and,
to this end, we see that a further modification of the integrand is necessary.
 Propositions (4.5) and (4.6) state  the conditions under which
 this modified integral is well defined and finite, respectively.

Note that if we define a function $ \rho_{H_{x,y}} $ on $G$ by
 $\rho_{H_{x,y}} := {\rho_{H^x \cap K^y} / \rho_{H^x}}
$, we have 
$$\rho_{H_{x,y}}(sz) = {\rho_{H^x \cap K^y}(sz)
 / {\rho_{H^x}(sz)}} = {{\Delta_{H^x \cap K^y} (s)}/ \Delta_{H^x} (s)}
\rho_{H_{x,y}}(z),$$
 for $s\in {H^x \cap K^y}$ and  $z\in G. $ Thus
$\rho_{H_{x,y}},$ restricted to $H^x ,$ is a $\rho$-function for the
homogeneous space $H^x/( {H^x \cap K^y}).$ We let $\mu_{H_{x,y}} $ be a 
quasi-invariant measure associated with this $\rho$-function and 
$ \lambda_{H_{x,y}} $ be the corresponding $\lambda $-function. 
Similarly, we can define a $\rho$-function $\rho_{K_{x,y}}$ for the 
homogeneous space $K^y/( {H^x \cap K^y})$ and the corresponding 
$\lambda $-function will be denoted by $ \lambda_{K_{x,y}}. $

\begin{pro} 
 Let $p,q$ and $m$ be positive real
numbers with $1\leq p,q'<\infty .$  
Then, for  $\sum_{i=1}^{\infty } f_i\otimes g_i \in L_p (\pi )\otimes^{\sigma }
L_{q'} (\gamma^*) $ and $x,y \in G,$  
\begin{eqnarray}
t\mapsto \sum_{i=1}^{\infty}{1\over {{\lambda_{H^x \cap K^y} (e,t)}^{1\over{m}}}}
 \lambda_H (x,t)^{1\over p}
 f_i (xt) {\otimes}_{x,y} \lambda_K (y,t)^{1\over q'} g_i (yt) 
\end{eqnarray}
\noindent  is a mapping on the coset space $G/({H^x \cap K^y })$
 in each of the following cases:
\begin{description}
 \item[(a)]  $p=m$ and $G/K$ having  invariant measure 
(or $q'=m$ and $G/H$ having  invariant measure); 
 \item[(b)]  $G/K$ and $G/H$  both having  invariant measures; 
\item[(c)]  $p=q'=m;$  
\item[(d)]  $H^x /( {H^x \cap K^y })$ and $K^y /( {H^x \cap K^y })$ 
having invariant measures. 
\end{description}
\end{pro}

\noindent
Proof: First let us consider the expression 
$${{{\lambda_H (x,t)}^{1\over p} {\lambda_K (y,t)}^{1\over {q'}}}
 \over {\lambda_{H^x \cap K^y} (e,t)}^{1\over m}} ,$$
under the cases (a), (b) and (c).

Consider (a). Assuming $p=m$ and using the identity (13), we have,
$$ \biggl({{\lambda_H (x,t)} \over {\lambda_{H^x \cap K^y} (e,t)
}}\biggr)^{1\over p}=
{\lambda_K (y,t)}^{-{1\over p}}.$$
If the measure on $G/K$ is invariant, then ${\lambda_K (y,t)} =1;$ hence
$${{\lambda_H (x,t)^{1\over p} {\lambda_K (y,t)}^{1\over {q'}}} \over 
{\lambda_{H^x \cap K^y} (e,t)}^{1\over p}}
 = \lambda_K (y,t)^{{1\over {q'}}- 
{1\over p}} =1.$$ A similar argument holds in the case where $q'= m $ and 
$G/H$ possesses an invariant measure.
In the case of (b), $\lambda_H (x,t) = \lambda_K (y,t) =1,$ and then 
by identity (13), $ \lambda_{H^x \cap K^y} (e,t) =1,$ giving 
$${{{\lambda_H (x,t)}^{1\over p} {\lambda_K (y,t)}^{1\over {q'}}}
 \over {\lambda_{H^x \cap K^y} (e,t)}^{1\over m}} =1.$$

Clearly, under condition (c),
$${{{\lambda_H (x,t)}^{1\over p} {\lambda_K (y,t)}^{1\over {q'}}}
 \over {\lambda_{H^x \cap K^y} (e,t)}^{1\over m}}=\biggl({{{\lambda_H (x,t)}
 {\lambda_K (y,t)}}
 \over {\lambda_{H^x \cap K^y} (e,t)}}\biggr)^{1\over p} =1,$$
using  the identity (13).

Therefore, under the conditions (a), (b) or (c), (30) can be simplified to
$$t\mapsto \Sigma_{i=1}^{\infty} f_i (xt) {\otimes}_{x,y}  g_i (yt) ,$$
 which is constant on each coset of 
${H^x \cap K^y }$ in $G .$ Hence it is a mapping on the coset space 
$G/({H^x \cap K^y }).$

 For the  case (d), it only remains to show that
\begin{eqnarray*}
{1\over {\lambda_{H^x \cap K^y} (e,st)}^{1\over m}} \lambda_H 
(x,st)^{1\over p} \lambda_K (y,st)^{1\over q'}&=&
{1\over {\lambda_{H^x \cap K^y} (e,t)}^{1\over m}} \lambda_H 
(x,t)^{1\over p} \lambda_K (y,t)^{1\over q'},
\end{eqnarray*}
for $s\in {H^x \cap K^y }.$  Letting $s = x^{-1}hx = y^{-1}ky,$ 
for $h\in H$ and $k\in K,$
we have
\noindent 
\begin{eqnarray}
\lefteqn{{1\over {\lambda_{H^x \cap K^y} (e,st)}^{1\over m}
} \lambda_H (x,st)^{1\over p}\lambda_K (y,st)^{1\over q'}}\nonumber \\
 &=&{\Bigl(}{\Delta_G (s)\over {\Delta_{H^x \cap K^y}(s)}}{\Bigl)}^{1\over m}
\Bigl({\Delta_H (h)\over {\Delta_G (h)}}{\Bigl)}^{1\over p}\Bigl({\Delta_K 
(k)\over {\Delta_G (k)}}{\Bigl)}^{1\over q'}
{1\over {{\lambda (e,t)}^{1\over m}}}\lambda_H (x,t)^{1\over p} \lambda_K (y,t)^
{1\over q'} .
\end{eqnarray}
 Since we 
assume that $H^x /({H^x \cap K^y })$ has  invariant measure,
 we have (see Reiter\cite{rei} p.159),
$$  \lambda_{H_{x,y}}(e,s) =  {{\rho_{H_{x,y}} (s)}\over 
 {\rho_{H_{x,y}} }(e)} = {{\Delta_{H^x \cap K^y} (s)}\over {\Delta_{H^x} (s)}}=
1.$$
Now $H$ and $H^x $ are closed conjugate subgroups of $G$ under an inner
automorphism $\tau : G\mapsto G$ given by $\tau (y) = x^{-1} y x.$
 Since $\tau $ is a topological 
isomorphism of $H$ onto $H^x $ we have $\Delta_H = \Delta_{H^x } \tau .$ 
This implies that $\Delta_{H^x} (h^x) = \Delta_H (h).$ Hence we have
\begin{eqnarray}
\quad \quad \quad \quad \quad \quad \quad {\Delta_H (h)\over \Delta_{H^x \cap K^y}  (s) }&=&{\Delta_K (k)\over 
\Delta_{H^x  \cap K^y}  (s)}=1,
\end{eqnarray}
 for $s\in H^x \cap K^y$ with $s = x^{-1}hx = y^{-1}ky.$
Considering the identity (17) and using the fact that $H^x /({H^x \cap K^y})$
and $K^y /({H^x \cap K^y})$ possess invariant measures, we have
\begin{eqnarray}
{\Bigl(}{\Delta_G (s)\over {\Delta_{H^x \cap K^y}(s)}}{\Bigl)}^{1\over m}
\Bigl({\Delta_H (h)\over {\Delta_G (h)}}{\Bigl)}^{1\over p}\Bigl({\Delta_K 
(k)\over {\Delta_G (k)}}{\Bigl)}^{1\over q'}&=&
{\Bigl(}{\Delta_G (s)\over {\Delta_{H^x \cap K^y}(s)}}{\Bigl)}^{{1\over m}-
{1\over p}-{1\over q'}}=1.
\end{eqnarray}
Thus, (31) simplifies to
\begin{eqnarray}
{{\lambda_H (x,st)^{1\over p}\lambda_K (y,st)^{1\over q'}}\over 
{{\lambda_{H^x  \cap K^y} (e,st)}^{1\over m}}} 
&=&{{\lambda_H (x,t)^{1\over p}\lambda_K (y,t)^{1\over q'}}\over 
{{\lambda_{H^x  \cap K^y} (e,t)}^{1\over m}}} ,
\end{eqnarray}
for $s\in {H^x  \cap K^y}$ and therefore,  (30) is a well 
defined  mapping in case (d) as well, completing the proof of the Proposition.
 
\begin{flushright}  \large$\diamond$\\
\end{flushright}

Recall, from the discussion preceding Lemma 2.2, that $\Upsilon $ denotes the set of all double cosets $H\times K :\Delta $
of $G\times G.$ For $x,y\in G,$ let $$M_{x,y}^{({{q'}\over {{q'}-1}})} =
 \int_{H^x \over  {H^x \cap K^y }} \lambda_{H_{x,y}} (e, \alpha ) 
d\mu_{H_{x,y}} (\alpha ) {\hbox{ and }} N_{x,y}^{({p\over {p-1}})} =
 \int_{K^y \over {H^x \cap K^y }} \lambda_{K_{x,y}} (e, \xi ) 
d\mu_{K_{x,y}} (\xi ).$$

\begin{pro}
For  $\sum_{i=1}^{\infty } f_i\otimes g_i \in L_p (\pi )\otimes^{\sigma }
L_{q'} (\gamma^*) $ the  integral 
\begin{eqnarray}
\int_{G\over {H^x \cap K^y}}\sum_{i=1}^{\infty}{1\over 
{{\lambda_{H^x \cap K^y} (e,t)}}}
 \lambda_H (x,t)^{1\over p}
 f_i (xt) {\otimes}_{x,y} \lambda_K (y,t)^{1\over q'} g_i (yt) d\mu_{x,y} (t) 
\end{eqnarray}
 is finite  for almost all $D(x,y)\in 
\Upsilon $ in each of the following cases:
\begin{description}
 \item[(a)]  $p=1$ and $G/K$ having finite invariant measure 
(or $q'=1$ and $G/H$ having finite invariant measure); 
 \item[(b)]  $G/K$ and $G/H$  both having finite invariant measures; 
\item[(c)]  $p=q'=1;$  
\item[(d)]  $1< p,q' <\infty $ with $1/p + 1/{q'} \geq 1$ and,
 $H^x /( {H^x \cap K^y })$ and $K^y /( {H^x \cap K^y })$ 
being compact for almost all $x,y \in G$ with $(x,y)\mapsto M_{x,y} N_{x,y}$
being a bounded function from $\Upsilon $ to ${\cal R}.$
\end{description} 
\end{pro}

\noindent
Proof: First let us consider the cases (a), (b) and (c).
 Using the disintegration of measures in the spaces 
involved (as discussed in the proof of Lemma 2.2 ), we get
 \begin{eqnarray*}
\lefteqn{\sum_{i=1}^{\infty} \int_{D(x,y)\in \Gamma}
\int_{G \over {H^x\cap K^y}} \| f_i(xt)\| 
\| g_i(yt)\|  d\mu_{D}(\underline{t})d\mu_{H,K} (D) }{\hspace{1.5in}} \\
 &=&\sum_{i=1}^{\infty} \int_{G\over H} 
\int_{G\over K}\| f_i (x)\| 
\| g_i (y))\| d\mu_H (x) d\mu_K (y),\\
&=&\sum_{i=1}^{\infty} \|f_i\|_1 \|g_i\|_1 .
\end{eqnarray*}
\noindent Now in the case of (a), (b) or (c), we know that 
$\sum_{i=1}^{\infty} \|f_i\|_1 \|g_i\|_1 \leq \sum_{i=1}^{\infty} \|f_i\|_p
 \|g_i\|_{q'} .$
Hence we have the desired result since $\sum_{i=1}^{\infty} f_i\otimes g_i \in  L_p (\pi )\otimes_G^{\sigma }
L_{q'} (\gamma^*).$

Now let us consider the case (d). In the remainder of the proof, 
${\lambda_{H^x \cap K^y} (\cdot ,\cdot )}$ will be written 
as ${\lambda (\cdot ,\cdot )},$ for  simplicity of notation. 
 Using the identity
(13) we see that 
\begin{eqnarray}
{{\lambda_H (x,t)^{1\over p} 
\lambda_K (y,t)^{1\over {q'}}}\over {\lambda (e,t)}} 
&=& \biggl( {{\lambda_K (y,t)}\over {\lambda (e,t)}}\biggr)^{1\over{ p'}}
\biggl( {{\lambda_H (x,t)}\over {\lambda (e,t)}}\biggr)^{1\over q} . 
\end{eqnarray}
Let $1/p + 1/q' -1 =1/r.$ Then ${1/{ p'}}= 1 - 1/p = 1/q' -1/r =1/q' (1-q'/r).$ Similarly,
$1/q= 1/p (1-p/r).$ Therefore,
\begin{eqnarray}
{{\lambda_H (x,t)^{1\over p} 
\lambda_K (y,t)^{1\over {q'}}}\over {\lambda (e,t)}} & =&
\biggl( {{\lambda_K (y,t)}\over {\lambda (e,t)}}\biggr)^{{1\over {q'}}
 (1-{{q'}\over r})}
\biggl( {{\lambda_H (x,t)}\over {\lambda (e,t)}}\biggr)^{{1\over p} (1-
{p\over r})} . 
\end{eqnarray}
Hence we have
\begin{eqnarray}
\lefteqn{I_i (x,y)=\int_{G\over {H^x \cap K^y}}
{{\lambda_H (x,t)^{1\over p} \lambda_K (y,t)^{1\over {q'}}}
\over {\lambda_{H^x \cap K^y} (e,t)}}
 \|  f_i (xt)\|   \| g_i (yt)\| d\mu _{x,y}(t) }{\hspace{0.5in}} \nonumber  \\
& =& \int_{G\over {H^x \cap K^y}} (\|  f_i (xt)\|^p   
\| g_i (yt)\|^{q'} )^{1\over r} \biggl(\biggl({{\lambda_H (x,t)}\over 
{\lambda (e,t)}}\biggr)^{1\over p} \|  f_i (xt)\|\biggr)^{1-{p\over r}} \times
\nonumber \\
& & \quad \quad \quad \quad \quad \quad \quad 
 \quad  \biggl(\biggl({{\lambda_K (y,t)}\over 
{\lambda (e,t)}}\biggr)^{{1\over {q'}}}
\| g_i (yt)\|\biggr)^{ 1-{{q'}\over r}}  d\mu _{x,y}(t) .
\end{eqnarray}
Using Corollary 12.5 of Hewitt and Ross\cite{hew}, the above can be simplified to
obtain
\begin{eqnarray}
\lefteqn{ I_i (x,y)  \leq  \biggl( \int \| 
 f_i (xt)\|^p   
\| g_i (yt)\|^{q'}  d\mu _{x,y}(t) \biggr)^{1\over r} \times }  \nonumber  \\
& &\biggl( \int {{\lambda_H (x,t)}\over 
{\lambda (e,t)}} \|  f_i (xt)\|^p d\mu _{x,y}(t) \biggr)^{{{q'}-1}
\over {q'}}  
 \biggl( \int  
{{\lambda_K (y,t)}\over {\lambda (e,t)}}
\| g_i (yt)\|^{q'} d\mu _{x,y}(t)\biggr)^{{p-1}\over p} .
\end{eqnarray}
\noindent where the three integrals are over the coset space 
${G/{H^x \cap K^y}}.$ Let us consider $\int_{G\over{H^x \cap K^y}}\Bigl ({1\over {\lambda
 (e,t)}} 
\lambda_H (x,t)\parallel f_i (xt)\parallel^p \Bigl ) d\mu _{(x,y)}(t).$ 
By Lemma 2.1, there exists a  quasi-invariant measure $\mu_{H_{x,y}} $ 
  on 
${H^x/({H^x \cap K^y})}$ such that
\begin{eqnarray*}
\lefteqn{\int_{G\over {H^x \cap K^y}}\Bigl ( 
{{\lambda_H (x,t)}\over {\lambda (e,t)}} \parallel
 f_i (xt)\parallel^p \Bigl ) d\mu _{(x,y)}(t) }\hspace{0.5in}  \\
& =&\int_{G\over {H^x}} \int_{H^x \over {H^x \cap K^y}}
\Bigl ( {\lambda (\alpha ,t)\over {\lambda_{H^x} (\alpha ,t)}}{\lambda_H (x,\alpha t)
\over {\lambda (e,\alpha t)}}  \parallel f_i (x\alpha t)\parallel^p \Bigl )
 d\mu_{H_{x,y}}(\alpha  )d\mu_{H^x}(t).
\end{eqnarray*}
For
$\alpha = x^{-1} hx$ with $h\in H,$ we get 
$${\lambda (\alpha ,t)\over {\lambda_{H^x} (\alpha ,t)}}{\lambda_H (x,\alpha t)
\over {\lambda (e,\alpha t)}}
={ \lambda_{H^x} (e,\alpha  )\over {\lambda (e,\alpha )}} =\lambda_{H_{x,y}}(e,\alpha ),$$\noindent 
\noindent $\lambda_{H_{x,y}}$ being a $\lambda $-function for 
${{H^x}/({H^x \cap K^y})}$ corresponding to the measure 
$\mu_{H_{x,y}}.$
Using the assumption that ${{H^x}/({H^x \cap K^y})}$ is compact
and the fact that  $\lambda_{H_{x,y}}(e,\alpha ) $ is bounded on compact 
sets (see property (c) on $\lambda$-functions,page 2), we have
\begin{eqnarray*}
\int_{H^x \over {H^x \cap K^y}}\lambda_{H_{x,y}}(e,\alpha )
 d\mu_{H_{x,y}}(\alpha  ) &=& M_{x,y}^{({{q'}\over {{q'}-1}})} < \infty .
\end{eqnarray*}
\par \noindent Thus
\begin{eqnarray}
\lefteqn{\int_{G\over{H^x \cap K^y}}\Bigl ( 
{{\lambda_H (x,t)}\over {\lambda (e,t)}}\|
 f_i (xt)\|^p \Bigl ) d\mu _{(x,y)}(t) }\hspace{.5in} \nonumber  \\
&\leq & M_{x,y}^{({{q'}\over {{q'}-1}})}\int_{G\over H^x}  \| f_i (xt)\|^p  d\mu_{H^x}(t) \nonumber \\
&= & M_{x,y}^{({{q'}\over {{q'}-1}})}\int_{G\over H} 
 \| f_i ( t)\|^p  d\mu_{H}(t) = M_{x,y}^{({{q'}\over {{q'}-1}})}\|f_i \|_p^p .
\end{eqnarray}
\par \noindent Similarly, if $K^y/{(H^x \cap K^y )}$ is compact,
\begin{eqnarray}
\int_{G\over {H^x \cap K^y}}\Bigl ({\lambda_K (y,t)\over 
{\lambda (e,t)}}\parallel  g_i (yt)\parallel^{{q'}} \Bigl ) d\mu _{(x,y)}(t)
&\leq & N_{x,y}^{({p\over {p-1}})}\|g_i \|_{q'}^{q'} . 
\end{eqnarray}
The inequalities (39), (40) and (41 ) imply that
\begin{eqnarray}
 I_i (x,y) &\leq & \biggl( \int_{G\over {H^x \cap K^y}}\|  f_i (xt)\|^p   
\| g_i (yt)\|^{q'}  d\mu _{x,y}(t) \biggr)^{1\over r}\times \nonumber \\
& & \quad \quad \quad \quad  M_{x,y}\|f\|^
{p({{{q'}-1}\over 
{q'}} )} N_{x,y}\|g\|^{{q'}({{p-1}\over p})}.
\end{eqnarray}
Note that
\begin{eqnarray*}
\lefteqn{\Biggl(\int_{D(x,y)\in \Upsilon }\biggl(\int_{G\over {H^x \cap K^y }}
\sum_{i=1}^{\infty} {{\lambda_H (x,t)^{1\over p} \lambda_K (y,t)^{1\over {q'}}}
\over {\lambda_{H^x \cap K^y} (e,t)}}
\|f_i(xt)\| \|g_i (yt) )\| d\mu_{x,y} (t) \biggr)^r
d\mu_{H,K} (D) \Biggr)^{1\over r} }\\
& = & \Biggl(\int_{ \Upsilon }\biggl(\sum_{i=1}^{\infty} 
\int_{G\over {H^x \cap K^y }} {{\lambda_H (x,t)^{1\over p} \lambda_K (y,t)^{1\over {q'}}}
\over {\lambda_{H^x \cap K^y} (e,t)}}
 \|f_i(xt) \| \| g_i (yt) )\| d\mu_{x,y} (t) \biggr)^r
d\mu_{H,K} (D) \Biggr)^{1\over r}\\
& = & \Biggl(\int_{D(x,y)\in \Upsilon }\biggl(\sum_{i=1}^{\infty} 
I_i (x,y) \biggr)^r
d\mu_{H,K} (D) \Biggr)^{1\over r} .
\end{eqnarray*}
Using generalised Minkowski's inequality (see Dunford and Schwartz\cite{dun},
p.529) we see that
\begin{eqnarray}
\Biggl(\int_{D\in \Upsilon }\biggl(\sum_{i=1}^{\infty} 
I_i (x,y) \biggr)^r
d\mu_{H,K} (D) \Biggr)^{1\over r} 
& \leq & \sum_{i=1}^{\infty} \Biggl(\int_{D\in \Upsilon }\biggl(I_i (x,y) \biggr)^r
d\mu_{H,K} (D) \Biggr)^{1\over r} 
\end{eqnarray}
\noindent Let ess sup${}_{D(x,y)}\{ (M_{x,y}N_{x,y})^r \} =S^r.$ Then, by 
(42) and (43) we have
\begin{eqnarray}
\lefteqn{\Biggl(\int_{D(x,y)\in \Upsilon }\biggl(\sum_{i=1}^{\infty} 
I_i (x,y) \biggr)^r
d\mu_{H,K} (D) \Biggr)^{1\over r} }   \nonumber    \\
& \leq &\sum_{i=1}^{\infty} \Biggl(\|f_i\|^{rp({{{q'}-1}\over {q'}})}
 \|g_i\|^{r{q'}({{p-1}\over p})} \times \nonumber \\
& &  \int_{D(x,y)\in \Upsilon } (M_{x,y} N_{x,y})^r
\biggl(\int_{G\over {H^x \cap K^y }}
 \|f_i(xt)\|^p \|g_i(yt)\|^{q'} d\mu_{x,y} (t)\biggr) d\mu_{H,K} (D) 
\Biggr)^{1\over r}, \nonumber \\
&\leq & \sum_{i=1}^{\infty}\Biggl( S^r \|f_i\|_p^{rp({{{q'}-1}\over {q'}})} \|g_i\|_{q'}^{r{q'}({{p-1}
\over p})} 
\int_{G\over H} \int_{G\over K}  \|f_i(x)\|^p \|g_i(y)\|^{q'}d\mu_H (x)
d\mu_K (y) \Biggr)^{1\over r} \\
&=& \sum_{i=1}^{\infty} \Biggl( S^r \|f_i\|_p^{p+{rp({{{q'}-1}\over {q'}})}} 
\|g_i\|_{q'}^{{q'}+{r{q'}({{p-1}
\over p})}}\Biggr)^{1\over r} .
\end{eqnarray}
where (44) is obtained using disintegration of measures (see proof of
 Lemma 2.2).
Since $p/r +p({q'}-1)/ {q'}= p(1/r +1 -1/{q'}) =p(1/p)=1,$ and similarly 
${q'}/r+{q'}(p-1)/ p =1,$ we obtain
\begin{eqnarray}
\Biggl(\int_{D(x,y)\in \Upsilon }\biggl(\sum_{i=1}^{\infty} 
I_i (x,y) \biggr)^r
d\mu_{H,K} (D) \Biggr)^{1\over r} \leq 
S\sum_{i=1}^{\infty}\|f_i\|_p \|g_i\|_{q'}. 
\end{eqnarray}
This proves the finiteness of the integral (35) for almost all
$D(x,y)\in \Upsilon $ under condition (d) together with
$1/p +1/{q'} >1.$

Now let us consider the case (d) 
 together with
 $1/p +1/{q'} =1.$ 
  \par \noindent Using H\"{o}lder's inequality, we get
\begin{eqnarray*}
\lefteqn{\int_{G\over {H^x \cap K^y}
}\sum_{i=1}^{\infty } {1\over {{\lambda (e,t)}}} \lambda_H (x,t)^{1\over p} 
\lambda_K (y,t)^{1\over {p'}}\|  f_i (xt)\|   \| g_i (yt)\|
d\mu _{x,y}(t) } \hspace{0.3in } \\
& \leq & \sum_{i=1}^{\infty } \Biggl 
(\int_{G\over{{H^x \cap K^y}}
}\Bigl ( {\lambda_H (x,t)^{1\over p}\over {{\lambda (e,t)}^
{1\over p}}} 
 \parallel f_i (xt)\parallel \Bigl )^p d\mu _{x,y}(t)\Biggl )^{1\over p} 
\times \\
& &\quad \quad \quad \quad \quad \quad \quad \quad \quad 
\Biggl (\int_{G\over 
{H^x \cap K^y}}
\Bigl ({\lambda_K (y,t)^{1\over {p'}}\over  
{{\lambda (e,t)}^{1\over {{p'}}}}}\parallel g_i (yt)\parallel \Bigl )^{{p'}}
d\mu _{x,y}(t)
\Biggl )^{1\over {{p'}}},\quad \quad \\
& = & \sum_{i=1}^{\infty } \Biggl (\int_{G\over {H^x \cap K^y}
}\Bigl ( {\lambda_H (x,t)\over {\lambda (e,t)}} \parallel 
 f_i (xt)\parallel^p \Bigl ) d\mu _{(x,y)}(t)\Biggl )^{1\over p} \times\\
& &\quad \quad \quad \quad \quad \quad \quad \quad \quad 
\Biggl ( \int_{G\over 
{H^x \cap K^y}}\Bigl ({\lambda_K (y,t)\over 
{\lambda (e,t)}}\parallel  g_i (yt)\parallel^{{p'}} \Bigl ) 
d\mu _{(x,y)}(t)\Biggl)^{1\over {{p'}}}.
\quad \quad
\end{eqnarray*}
By (40) and (41) we have \par \noindent 
\begin{eqnarray*}
\lefteqn{\int_{G\over {H^x \cap K^y}
}\sum_{i=1}^{\infty } {{\lambda_H (x,t)^{1\over p} 
\lambda_K (y,t)^{1\over {p'}}}\over {\lambda (e,t)}} \parallel  f_i (xt)\parallel   \parallel  g_i (yt)\parallel 
d\mu _{(x,y)}(t) }\\
&\leq & M_{x,y} N_{x,y}\sum_{i=1}^{\infty } 
\|f_i\|_p\| g_i \|_{p'}.
\end{eqnarray*}
 
Again, since $\sum_{i=1}^{\infty} f_i\otimes g_i \in  L_p (\pi )\otimes_G^{\sigma }
L_{q'} (\gamma^*),$ 
we see that $\sum_{i=1}^{\infty} \|f\|_p \|g\|_{q'}
<\infty .$
 Hence the result follows.

\begin{flushright}  \large$\diamond$\\
\end{flushright}

In view of Propositions 4.5 and 4.6, we can formally define 
the convolution of functions in $L_p(\pi )$ and $L_{q'} (\gamma^* ).$ 

\begin{dfn} 
Let $H$ and $K$ be regularly related. For each  $x,y \in G$ let $\mu _{x,y}$ be a quasi-invariant 
measure  on the homogeneous space ${G/{(H^x \cap K^y )}}$
 so that the identity (9) holds. Let $p,q$  be positive real numbers.  
The map $\Psi$ on 
$L_p(\pi ) {\otimes_{\sigma }} L_{q'} (\gamma^* )$ is defined by  
\begin{eqnarray}
\lefteqn{(\Psi (\Sigma_{i=1}^{\infty}f_i\otimes g_i))(x,y) :=}     \\
&  &\int_{G\over {H^x \cap K^y}
}\Sigma_{i=1}^{\infty}{1\over {{\lambda_{H^x \cap K^y} (e,t)}}}
 \lambda_H (x,t)^{1\over p}
 f_i (xt) {\otimes}_{x,y} \lambda_K (y,t)^{1\over q'} g_i (yt)
 d\mu _{(x,y)}(t) \nonumber
\end{eqnarray}
\noindent  whenever one of the following conditions holds:
\begin{description}
 \item[(a)]  $p=1$ and $G/K$ has finite invariant measure 
(or $q'=1$ and $G/H$ has finite invariant measure); 
 \item[(b)]  $G/K$ and $G/H$  both have finite invariant measures; 
\item[(c)]  $p=q'=1;$  
\item[(d)]  $1< p,q' <\infty $ with $1/p + 1/{q'} \geq 1,$ $H^x /( {H^x \cap K^y })$ and $K^y /( {H^x \cap K^y })$ 
are compact and possess invariant measures for almost all $x,y \in G$ and the map $(x,y)\mapsto M_{x,y}
N_{x,y} $ is bounded from $\Upsilon $ to ${\cal R}.$
\end{description}
\end{dfn}

It is clear that for each $(x,y)\in G\times G,$ the value of 
${(\Psi (\Sigma_{i=1}^{\infty}f_i\otimes g_i))(x,y) }$ belongs to the quotient
space ${\cal A}_{x,y} .$ We investigate the properties of the image space of  
$\Psi$ in the following section.

\subsection{The space $A_p^q$}
First we shall show that the image space of $\Psi $ 
consists of mappings which are constant on the right cosets $(G\times G)
/{\Delta} $
under the conditions (a), (b) and (c) of Definition 4.7.
\noindent
\begin{pro}Let $\alpha $ be an element of the image space of $\Psi .$
 For any \newline \noindent $h_0 \in H, k_0 \in K, $
  $x,y \in G$ and  $s\in G/{(H^x \cap K^y )} $
$$\alpha (h_0 xs,k_0 ys)=\pi (h_0 )\otimes 
\gamma^* (k_0 ) \alpha (x, y)$$
\noindent under the conditions (a), (b) and (c) of Definition 4.7. 

\end{pro}

\noindent 
Proof: By Proposition 4.4 we have 
$${\cal H}_{xs,ys} =  {\cal H}_{x,y}$$  for all $x,y,s\in G.$

Now any element $\omega \otimes_{x,y} \varrho $ of ${\cal A}_{x,y}$ is of 
the form
\noindent
\begin{eqnarray*}
\lefteqn{\omega \otimes_{x,y} \varrho = {\cal H}_{x,y} + \omega \otimes 
 \varrho }\hspace{0.2in} \\
&=&\langle \{ {\pi }^{x} (b) \xi \otimes \eta - \xi \otimes 
({\gamma}^{y }(b))^* 
\eta ,b \in H^{x} \cap K^{y} , \xi \in H(\pi ), \eta \in H(\gamma )\}\rangle  + 
\omega \otimes \varrho .
\end{eqnarray*}
\noindent
If this element is translated by $\pi (h_0 )\otimes \gamma^* (k_0 ) $ from 
the left,
we get  \medskip \par \noindent
{\leftline{$ \pi (h_0 )\otimes \gamma^* (k_0 )(\omega \otimes_{x,y} \varrho 
) $}}
 $$  =\langle \{ \pi (h_0 )
{\pi }^{x} (b) \xi \otimes \gamma^* (k_0 )\eta -\pi (h_0 ) \xi \otimes 
\gamma^* (k_0 )
({\gamma}^{y }(b))^* \eta  ;\  b \in H^{ x} \cap K^{ y}
\}\rangle +\pi (h_0 ) \omega \otimes \gamma^* (k_0 )  \varrho .$$
But $$\pi (h_0 ){\pi }^{x} (b) \xi =\pi^ {{h_0 }x}(b){\pi }(h_0 ) \xi 
{\hbox
 { and }}
 \gamma^* (k_0 )({\gamma}^{y }(b))^* \eta = ({\gamma}^{k_0 y }(b))^* \gamma^* 
(k_0 )\eta ,$$
hence \par \noindent 
\begin{eqnarray}
\lefteqn{\pi (h_0 )\otimes \gamma^* (k_0 )(\omega \otimes_{x,y} \varrho )
  }\nonumber \\ 
& = &\langle \{ \pi^ {{h_0 }x}(b){\pi }(h_0 ) \xi \otimes \gamma^* 
(k_0 )\eta -\pi (h_0 ) \xi \otimes ({\gamma}^{k_0 y }(b))^* \gamma^* 
(k_0 )\eta 
;\ b \in H^{h_0 x} \cap K^{k_0 y} \}\rangle \nonumber \\
& & \quad  \quad  \quad  \quad  \quad  \quad  \quad  \quad  \quad  
+\pi (h_0 ) \omega \otimes \gamma^* (k_0 )  \varrho ,\nonumber \\
&= &{\cal H}_{ h_0 x,k_0 y}  + \pi (h_0 ) \omega  \otimes \gamma^* (k_0 ) 
 \varrho  ,\nonumber \\
& = & \pi (h_0 ) \omega  \otimes_{ h_0 x,k_0 y} \gamma^* (k_0 )  \varrho .
\end{eqnarray}

 Any $\alpha $ in the image space of $\Psi $  can be expressed as
$\Psi (\sum_{i=1}^{\infty } f_i \otimes g_i ).$ Without loss of generality,
we consider  an element of the form $\Psi ( f \otimes g );$ 
the argument is then valid for any $\alpha $ by linearity.
 Consider the
 homeomorphism $\phi_{s} :G/({H^x\cap K^y}) \mapsto G/({H^x\cap K^y})^s$
given by $\phi_{s} (v) = s^{-1} v,$  and use the fact that
$\mu_{xs,ys} =\mu_{x,y}\circ \phi_{s}$ to get
\begin{eqnarray*}
\lefteqn{(\Psi ( f \otimes g ))( h_0 xs,k_0 ys) }\hspace{0.1in} \\
& =&\int_{G\over {H^{xs} \cap K^{ys}}}
{1\over {\lambda_{(H^x \cap  K^y)^s} (e,t)}} 
\lambda_H ( xs,t)^{1\over p} f (xst) {\otimes}_{ xs,ys} 
\lambda_K ( ys,t)^{1\over {q'}}
 g (yst) d\mu _{(xs,ys)}(t) \\
& =&\int_{G\over {H^{xs} \cap K^{ys}}}
({{\lambda_{H^x \cap  K^y} (e,s)}\over 
{\lambda_{H^x \cap  K^y} (e,st)}})({ \lambda_H ( x,st)\over
 {\lambda_H ( x,s)}})^{1\over p} f (xst) {\otimes}_{ xs,ys} \\
& &\quad \quad \quad \quad \quad \quad \quad \quad \quad \quad \quad 
\quad \quad \quad \quad 
({ \lambda_K ( y,st)\over {\lambda_K ( y,s)}})^{1\over {q'}}
 g (yst) d\mu _{(xs,ys)}(t) ,\\
& =&\int_{G\over {H^{x} \cap K^{y}}}
({{\lambda_{H^x \cap  K^y} (e,s)}\over 
{\lambda_{H^x \cap  K^y} (e,t)}})
({ \lambda_H ( x,t)\over {\lambda_H ( x,s)}})^{1\over p} f (xt) {\otimes}_{ xs,ys} 
({ \lambda_K ( y,t)\over {\lambda_K ( y,s)}})^{1\over {q'}}
 g (yt) d\mu _{(x,y)}(t) ,\\
&=&{{{\lambda_{H^x \cap  K^y} (e,s)}}\over 
{{\lambda_H ( x,s)}^{1\over p} {\lambda_K ( y,s)}^{1\over {q'}}}}
 (\pi (h_0 ) \otimes \gamma^* (k_0 ))[\Psi (f
\otimes g)]( x, y) . 
\end{eqnarray*}
\noindent
  Hence we see that 
$$\alpha (h_0 xs,k_0 ys)=\pi (h_0 )\otimes 
\gamma^* (k_0 ) \alpha (x, y)$$
\noindent only if $${{{\lambda_{H^x \cap  K^y} (e,s)}}\over 
{{\lambda_H ( x,s)}^{1\over p} {\lambda_K ( y,s)}^{1\over {q'}}}}=1$$
\noindent for all $s\in G/{(H^x \cap K^y )}.$
It is clear that this last condition is true in the cases (a), (b) and (c)
given in Definition 4.7. 

\begin{flushright}  \large$\diamond$\\
\end{flushright}

\subsubsection*{The structure of the image space of $\Psi$}

The image space of $\Psi $ is contained in a space of mappings acting on 
$G\times G,$ whose values at $(x,y)\in G\times G$ belong to a 
collection  of Banach
 spaces \newline \noindent $\{{\cal A}_{x,y} : (x,y)\in {G\times G}\}.$
This suggests that the image space has the structure of the space of 
cross-sections of 
a Banach bundle or a Banach semi-bundle where the bundle space is a union of quotient 
spaces of a given Banach space. 

\noindent Let 
\begin{eqnarray*}
{\cal B}_0 &= &{\cal H}(\pi )\otimes^{\sigma }{\cal H}(\gamma^*)
 \times G\times G,\\
{\cal B}_0^{\Delta }&= & {\cal H}(\pi )\otimes^{\sigma }{\cal H}(\gamma^*) 
\times {({G\times G})/ \Delta },\\ 
{\cal B}_1 &= & \cup_{(x,y)\in {G\times G}}\{ {\cal H}_{x,y} \times 
\{(x,y)\}\},\\
 {\cal B}_1^{\Delta } &= & \cup_{(x,y)\in G\times G}\{ {\cal H}_{x,y} \times
\{(x,y) \Delta \}\},\\ 
{\cal B}_2 &= & \cup_{(x,y)\in {G\times G}}\{ {\cal A}_{x,y} \times 
\{(x,y)\}\},{\hbox{ and }}\\ 
{\cal B}_2^{\Delta }&= & \cup_{(x,y)\in {G\times G}}\{ {\cal A}_{x,y} \times
\{(x,y) \Delta \}\}.
\end{eqnarray*}
It is clear that ${\cal B}_1 $ is a subspace of ${\cal B}_0 ,$ 
and $ {\cal B}_1^{\Delta } $ is a subspace of ${\cal B}_0^{\Delta} .$

For $(x,y) \in G\times G,$ let $r(x,y) \in (G\times G)/{\Delta }$ be the right coset 
to which $(x,y)$ belongs.
  With $j$ denoting any one of
$\{0,1,2\},$ let 
$ \theta_{j} :{\cal B}_{j}\mapsto  {G\times G}$ be defined by
$\theta_{j} (\zeta ,(x,y))=(x,y),$ and let
$ \theta_{j}^{\Delta } :{\cal B}_{j}^{\Delta }\mapsto 
{({G\times G})/ {\Delta }}$ be defined by
$\theta_{j}^{\Delta } (\zeta ,(x,y){\Delta })=(x,y){\Delta },$ 
where $\zeta $ belongs to the corresponding Banach space. 
Let $q:{\cal B}_0 \mapsto {\cal B}_2$ be the quotient
 map given by $q(h,x) = (\{{\cal H}_x +h\},x).$ Similarly, the quotient
 map $q_{\Delta }:{\cal B}_0^{\Delta } \mapsto {\cal B}_2^{\Delta }$ is given by 
 $q_{\Delta }(h,r(x)) = (\{{\cal H}_x +h\},r(x)).$
${\cal B}_0 $ has the product topology, and we topologize 
 ${\cal B}_2^{\Delta }$ so that the map $p_{\Delta }$ is continuous and open.

Define $ {\underline{\cal B}}_{j} :=({\cal B}_{j},{\theta }_{j})$ and 
$ {\underline{\cal B}}_{j}^{\Delta } :=({\cal B}_{j}^{\Delta },{\theta }_{j}
^{\Delta }),$ for $j\in \{0,1,2\}.$

 The space ${{(G\times G)}/ \Delta}$ is Hausdorff since $\Delta $ is a 
closed subgroup of  ${G\times G}.$
It is can be easily seen that ${\underline{\cal B}}_{j}$ and 
${\underline{\cal B}}_{j}^{\Delta },$ 
for $j\in \{0,2\},$
are bundles over $G\times G$ and ${{G\times G}\over {\Delta }},$ respectively.
 Moreover,
 $ {\underline{\cal B}}_{0}$ and $ {\underline{\cal B}}_{0}^
{\Delta }$ are trivial bundles with constant fiber \newline 
\noindent ${\cal H}(\pi )
\otimes^{\sigma }{\cal H}(\gamma^*).$ Although ${\underline{\cal B}}_{2}$
 and ${\underline{\cal B}}_{2}^{\Delta }$ fail to be Banach bundles in general,
we see that they possess necessary properties
to become Banach semi-bundles. 
  For each $ z={\Delta }(x,y)\in  (G\times G)/{\Delta } , $ the fibre of  
  ${\cal B}_{2}^{\Delta } $ 
 over $z$ is 
${\cal B}_{2,z}^{\Delta } = \{{\cal A}_{x,y}\times \{z \}\}$.
We see that  ${\cal B}_{2,z}^{\Delta } $, 
is a Banach space with the norm $\| (\eta  , z) \|_{{\cal B}_{2,z}^{\Delta } } $
  defined by $\| (\eta  ,z)\|_{{\cal B}_{2,z}^{\Delta } } =
\| \eta  \| $ where  
$\| \eta  \| $ means the norm in  ${\cal A}_{x,y} .$   The
 operations $+$ and $.$ in  ${\cal B}_{2,z}^{\Delta } $ are defined, in an obvious manner, using
 $+$ and $.$ in  ${\cal A}_{x,y} .$
We can define and topologize the fibres ${\cal B}_{2,z} $ in ${\cal B}_{2} $
and define the operations $+$ and $.$  
in a  similar manner.

\begin{lem}
$ (\eta  , z)\mapsto   \| (\eta  , z)
 \|_{{\cal B}_{2,z}^{\Delta } } $  is upper semi-continuous on  
  ${\cal B}_{2}^{{\Delta }} $ to ${\cal R} .$ 
A similar result holds in the case of  ${\cal B}_{2}.$
\end{lem}

\noindent Proof: Let $\{ (\eta_i , z_i): i\in I \}$ 
be a net of elements in 
${\cal B}_{2}^{\Delta }$ with
$ (\eta_i , z_i)\rightarrow (\eta  , z ).$ Then there exist a 
sequence $\{ (\varphi_i , u_i) \} $ and an element $(\varphi ,u)$ in $
{\cal B}_0^{\Delta }$ such that \newline 
\noindent $q_{\Delta }((\varphi_i , u_i)) = (\eta_i , z_i)$ for all
$i\in I$ , $q_{\Delta }((\varphi , u)) = (\eta  , z)$ and $ (\varphi_i , u_i)
\rightarrow (\varphi , u ).$  Now since $
\| (\eta  , z)\|_{{\cal B}_z }=\| \eta  \|= 
inf_{h\in {\cal H}_{x,y}} \|\varphi + h \|,$ without loss of generality 
we can choose $\varphi $ such that, for a given $\epsilon >0,$ we have 
\begin{eqnarray}
\|\varphi \| &<& \|\eta  \|+ \epsilon .
\end{eqnarray}
Also,\begin{eqnarray}
\|\eta_i \| &\leq &\|\varphi_i \|
\end{eqnarray}
 for all $i\in I.$
Since $\|\varphi_i \|\rightarrow \|\varphi \|,$  then 
from (49) and (50) we have
$$\|\eta_i \|\leq \|\eta  \|+\epsilon ,$$
 for $i$ sufficiently large.

The proof is similar in the case of ${\cal B}_{2} .$

\begin{flushright} \large$\diamond$\\
\end{flushright}

\begin{lem} 
The operation $+$ is continuous on ${\cal B}_{2,z}^{\Delta } \times {\cal B}_{2,z}^{\Delta }
 $ to ${\cal B}_{2,z}^{\Delta } ,$
and  for each $\lambda $ in $C,$ the map $b\mapsto \lambda b$ is 
continuous on $ 
{\cal B}_{2}^{\Delta } $ to ${\cal B}_{2}^{\Delta }$ .
A similar results hold in  ${\cal B}_{2} .$ 
\end{lem}

\noindent Proof: Since the topology induced from ${\cal B}_{2}^{\Delta }$ on its fibres
is just the Banach space topology, the operations + and . are continuous.

\noindent Similarly for ${\cal B}_{2} .$

\begin{flushright}  \large$\diamond$\\
\end{flushright}

\begin{lem}          
If $z\in (G\times G)/{\Delta }$ and $\{b_i :i\in I\},$ is any net of elements in 
${\cal B}_{2}^{\Delta } $ 
such that 
 $\| b_i \|\rightarrow  0 $ and $\theta_2^{\Delta } (b_i ) \rightarrow z$,
 then $ b_i \rightarrow  0_z$ where $0_z $
 is the zero element in ${\cal B}_{2,z}^{\Delta } .$ A similar result holds in
 ${\cal B}_{2} .$
\end{lem}

\noindent Proof: Any element $b_i \in {\cal B}_{2}^{\Delta } $ is of the form
$b_i = (\omega_i + {\cal H}_{x_i,y_i } , (x_i ,y_i ){\Delta })$ where $x_i ,y_i \in G.$
 Since $ \| b_i \| = \inf \{\| \omega_i + h \|  : h\in {\cal H}_{x_i ,y_i} \},$ with
$ \omega_i \in {\cal H},$ there exists an $h_i  \in {\cal H}_{x_i ,y_i}
$ such that  $$\| \omega_i + h_i\|  <  \| b_i \| +1/2^i$$ for all $i\in I.$
This implies that the net of elements $\omega_i + h_i $ in  
${\cal H}(\pi )\otimes^{\sigma }{\cal H}(\gamma^*) $
 has the property that  $\omega_i + h_i  \rightarrow {\underline{0 }},$
${\underline{0 }}$ being the zero element in ${\cal H}(\pi )\otimes^{\sigma }
{\cal H}(\gamma^*) .$ If  $\theta_2^{\Delta }
 (b_i ) = (x_i ,y_i ){\Delta }
\rightarrow z ,$ this means that $b_i \rightarrow 0_z.$ 

\begin{flushright}  \large$\diamond$\\
\end{flushright}

\begin{lem}${\underline{\cal B}}_{2}^{\Delta } $ 
and ${\underline{\cal B}}_{2}$ 
are  Banach semi-bundles over ${(G\times G)/{\Delta }}$ and $(G\times G)$ respectively.
\end{lem}

\noindent Proof: 
The result follows from  Lemmas 4.9, 4.10 and 4.11.

\begin{flushright} \large$\diamond$\\
\end{flushright}

\begin{pro}
The Banach semi-bundle retraction $${\underline {\cal B}_2^{\Delta }}^{\# }=
({{\cal B}_2^{\Delta }}^{\# }, 
{\theta_2 ^{\Delta }}^{\# })$$ of ${\underline {\cal B}_2^{\Delta }}$ by $r$ is 
topologically 
equivalent to ${\underline {\cal B}_2}.$
\end{pro}

\noindent Proof:
Consider the diagram 
$${{\cal B}_0^{\Delta }}^{\# }\stackrel{q_{\Delta }^{\# } }{\mapsto }
{{\cal B}_{2}^{\Delta }}^{\# } \quad \quad \quad \quad \quad \quad \quad \quad $$
$$\imath \uparrow \quad \quad \quad \quad  \uparrow \jmath 
\quad \quad \quad \quad \quad \quad \quad \quad $$
$${{\cal B}_0}\stackrel{ q }{\mapsto }{\cal B}_{2} \quad \quad \quad \quad 
{{\cal B}_0^{\Delta }}\stackrel{ q_{\Delta } }{\mapsto }{\cal B}_{2}^{\Delta }$$
$$ \theta_0 \searrow  \swarrow  \theta_2  \quad \quad \quad \quad 
\theta_0^{\Delta } \searrow  \swarrow  \theta_2^{\Delta } $$
$$ G\times G  \quad \quad \stackrel{r}{\longrightarrow }\quad \quad (G\times G)/{\Delta } $$
where $q_{\Delta }^{\# } = j^{\# } (q_{\Delta } )$ (see Sec.1.2, p.4) and $\jmath  $ 
is defined so that $  q_{\Delta }^{\# } \circ \imath = \jmath 
\circ q , $ $\imath $ being the homeomorphism  stated in 
Proposition 1.3. It is clear that $q_{\Delta }^{\# }$ is the quotient map.
Hence, (also  by the discussion on pages 4 and 5,) 
$q_{\Delta }^{\# } $ is continuous and open. Obviously, $\jmath  $ defines a bijection from 
${{\cal B}_{2}}$ onto ${{\cal B}_{2}^{\Delta }}^{\# }.$  We need to show that $\jmath  $
and its inverse are continuous.
 Now $q_{\Delta } $ is open  by the definition of the topology of ${{\cal B}_{2}^{\Delta }}
$ and the right hand side of the above diagram commutes.
Since the maps $\imath $, $q_{\Delta }^{\# }$ and $q$ are continuous and open it is 
clear that $\jmath $ is continuous and open, as required. 

\begin{flushright}  \large$\diamond$\\
\end{flushright}

\begin{pro}
Let $f:G\times G\mapsto {\cal B}_2 $ be a continuous cross-section 
 which is constant on 
equivalence classes. Then the function $g$ defined by $g(p(x,y)) = 
{i^{\# }}(f (x,y)),$ where $(x,y)\in G\times G,$ is a 
continuous cross-section from $(G\times G)/{\Delta }$ to ${\cal B}_2^{\Delta }.$
\end{pro}

\noindent Proof:By Proposition 4.13, a continuous cross-section  $f$ 
of ${{\underline{\cal B}}_2 }$
 can be regarded as a cross-section of ${{\underline{\cal B}}_2^{\Delta }}^{\# }.$ 
 Define $g':(G\times G)\mapsto {{\cal B}_2^{\Delta }}$ so that
     $$ g'(x,y) = {i^{\# }}(f (x,y)) {\hbox { for all }} (x,y)\in G\times G.$$   
(See  the diagram below.)
$${{\cal B}_2 }={{\cal B}_2^{\Delta }}^{\# }\stackrel{i^{\# }}{\longmapsto }
{{\cal B}_2^{\Delta }}$$
$$f \nwarrow  \quad \nearrow g' $$
$$X$$
\noindent 

Consider the function $g:(G\times G)/{\Delta }\mapsto
{{\cal B}_2^{\Delta }}$ which factors through the diagram
 $$\quad \quad  {{\cal B}_2^{\Delta }}$$
$$g'\nearrow \quad \uparrow g$$
$$G\times G\stackrel{r}{\longmapsto }( G\times G)/{\Delta }$$
\noindent It is clear that  $g$ is well defined since $f$ is 
constant on the equivalence classes. Also,
 $g(p(x,y)) = g'(x,y)$ for any $(x,y)\in  G\times G$ and we see that $g(z)\in 
{\cal B}_{2,z}^{\Delta }$ for any $z\in ( G\times G)/{\Delta }.$ Hence $g$ is a cross-section of
${{\underline{\cal B}}_2^{\Delta }}.$ Moreover, it is continuous since $p$ is open.

\begin{flushright} \large$\diamond$\\
\end{flushright}

\begin{lem}
Consider the conditions (a), (b), (c) and (d) of  Definition
4.7.
For $\sum_{i=1}^{\infty}f_i\otimes g_i \in L_p (\pi )\otimes L_{q'} (
\gamma ),$ the element $\Psi (\sum_{i=1}^{\infty}f_i\otimes g_i)$ is a 
cross-section 
of  ${\underline {\cal B}}_2^{\Delta },$ if the integral (47)
is constructed under one of the conditions (a), (b) or (c). It is a cross-section 
of  ${\underline {\cal B}}_2 $ if it is constructed under the condition (d).
\end{lem}

\noindent Proof: This is an immediate consequence of Proposition 4.8.

\begin{flushright}  \large$\diamond$\\
\end{flushright}
\noindent
\begin{dfn} The \sf space $A_p^q$ \it  is defined to be the range of $\Psi $ with the
 quotient norm.
\end{dfn}

In other words, $A_p^q$ is contained in the space of cross-sections of the Banach semi-bundle
 ${\underline {\cal B}}_2^{\Delta }$ in the cases (a), (b) and (c) of 
Definition 4.7. In the case (d), it is contained in the space of 
cross-sections of the 
Banach semi-bundle  ${\underline {\cal B}}_2 .$

By a \sf continuous family of  functions \rm we mean a family of 
 functions $\{\beta_{x} : x\in G \}$ such that 
$(x,t)\mapsto \beta_{x} (t)$ is a continuous map from $G\times G$ to ${\cal R}.$
\begin{pro} Suppose that the spaces $G/H, G/K$ and the numbers $p,q$ 
 satisfy one of the conditions (a), (b), (c) or (d) as described in Definition
4.7. Suppose  further that 
 there exists a continuous family $\{\beta_{x,y} : 
(x,y)\in G\times G \}$ of 
 functions where $\beta_{x,y}$ is a Bruhat function for ${H^x \cap K^y } .$  
Let $f$ and $g$ be functions with compact support from $ L_p (\pi )$ and 
$L_{q'} (\gamma^* )$ respectively.
 Then, 
\begin{eqnarray}
(x,y)&\mapsto &\int_{G\over {H^x \cap K^y }}{1\over {\lambda_{H^x \cap  K^y} 
(e,t)}} \lambda_H ( x,t)^{1\over p} f(xt)\otimes_{x,y} \lambda_K
 ( x,t)^{1\over {q'}}g(yt)
d\mu_(x,y) (t) \nonumber \\
&  & \quad \quad       
\end{eqnarray}
is a continuous cross-section of the corresponding Banach semi-bundle.
\end{pro}

 \noindent 
Proof: It can be easily seen that
for any $x\in G$ and  $f\in L_p (\pi ),$ the function ${}_x f $ defined by 
${}_x f (t) = f(xt )$ is a function in $L_p (\pi^x ).$  Similarly, a function
$g\in L_{q'} (\gamma^* )$ gives rise to a function ${}_y g$ in $L_{q'} (\gamma^* ).$ \par 
Now suppose $f $ and $g$ are continuous with compact support. Then 
there exist compact sets $G_1$ and $G_2$ of $G$ such that $ H^x G_1$ and 
$ K^y G_2 $ are the supports of ${}_xf $ and ${}_y g$ respectively.
\par \noindent  
\noindent Suppose that the integral in (51) is constructed under one of the 
conditions (a), (b), (c) or (d) of Definition 4.7. 
 Consider the map 
$$(x,y,t)\mapsto \beta_{x,y} (t) 
 \lambda_H ( x,t)^{1\over {p}} 
f(xt)\otimes \lambda_K (y,t)^{1\over {q'}} g(yt)$$ from 
$(G\times G\times G )$ to ${\cal B}_0 .$ This is a  cross-section 
of the Banach bundle retraction of ${\underline{\cal B}}_0$ by $p:G\times 
G\times G\mapsto G\times G.$ It is a continuous cross-section since, under the 
assumptions,  $\{\beta_{x,y} : 
(x,y)\in G\times G \}$ is a continuous family of Bruhat functions.
Therefore, we can form the integral $${\widetilde {\Gamma }} (x,y) :=
\int_G \beta_{x,y} (t) 
 \lambda_H ( x,t)^{1\over p} 
f(xt)\otimes \lambda_K (y,t)^{1\over {q'}}g(yt) d\nu_G (t)$$  
  and by Lemma 1.4, we see that  
${\widetilde {\Gamma }}$  is a
continuous cross-section of ${\underline {\cal B}}_0 .$\\
Considering the diagram
$$ {\cal B}_0 \stackrel{q}{\mapsto }{\cal B}_2 $$
 $${\widetilde {\Gamma }} \uparrow \quad \quad \quad \nearrow \Gamma , $$
 $$G\times G \quad \quad $$
where $q(\xi ,x) = (\{ {\cal H}_{x} + \xi \} ,x ),$ we find that 
 $$\Gamma (x,y) := \int_G \beta_{x,y} (t) 
 \lambda_H ( x,t)^{1\over q}f(xt)\otimes_{x,y} 
\lambda_K (y,t)^{1\over {p'}} g(yt) d\nu_G (t) $$ is a
continuous cross-section of ${\underline {\cal B}}_2 .$ Note that the
property (i) of $\lambda$-functions on page 2 implies that we can assume ${\rho_{H^x \cap  K^y}(e)} =1,$
for all $x,y\in G.$
 Using Corollary 1.2, we get
\begin{eqnarray}                                       
\Gamma (x,y)
& = &\int_{G\over {H^x \cap K^y }} \int_{H^x \cap K^y } 
\beta_{x,y} (st) {1\over {{\lambda_{H^x \cap  K^y}(e,st)}}} 
\lambda_H ( x,st)^{1\over p} f(xst)\otimes_{x,y} \nonumber \\
& & \quad \quad  \quad \quad  \quad \quad \quad \quad  \quad \quad 
\lambda_K (y,st)^{1\over {q'}}g(yst)d\nu_H (s) d\mu_{x,y} (t).
\end{eqnarray} 
But under the conditions (a), (b), (c)
 or (d) in Definition 4.7, (34) implies that, for $s\in {H^x \cap  K^y},$ 
\begin{eqnarray*}
{{\lambda_H ( x,st)^{1\over {p}} \lambda_K (y,st)^{1\over {q'}}}\over 
{\lambda_{H^x \cap  K^y}(e,st)}}& =& {{\lambda_H ( x,t)^{1\over {p}} 
\lambda_K (y,t)^{1\over {q'}}}\over 
{\lambda_{H^x \cap  K^y}(e,t)}}.
\end{eqnarray*}
\noindent
Therefore the integral (52) can be simplified to give
\begin{eqnarray*}
\Gamma (x,y)
& = & \int_{G\over {H^x \cap K^y }}{1\over {{\lambda_{H^x \cap  K^y} 
(e,t)}}} \lambda_H ( x,t)^{1\over p} f(xt)\otimes_{x,y} \lambda_K
 ( x,t)^{1\over {q'}}g(yt)\times \\
& & \quad \quad \quad \quad \quad \quad \quad \quad \quad \quad 
\quad  \Biggl( \int_{H^x \cap K^y }
\beta_{x,y} (st) d\nu_{H^x \cap  K^y} (s) \Biggr)d\mu_{x,y} (t),\\
&=&\int_{G\over {H^x \cap K^y }}
{1\over {\lambda_{H^x \cap  K^y} 
(e,t)}} \lambda_H ( x,t)^{1\over p} f(xt)\otimes_{x,y} \lambda_K
 ( x,t)^{1\over {q'}}g(yt)d\mu_{x,y} (t).
\end{eqnarray*}

\noindent
Hence the mapping given by (51) is continuous in the Banach semi-bundle
 ${\underline{\cal B}}_2 .$
In the case of (a), (b) or (c) in Definition 4.7, we can consider the  
mapping (51)
as a cross-section of the Banach semi-bundle retraction
 ${{\underline{\cal B}}_2^{\Delta}}^{\# }$ of
 ${\underline{\cal B}}_2^{\Delta}$ by the cannonical mapping $r:G\times G
\rightarrow {({G\times G})/ {\Delta}} .$
By Proposition 4.14, this cross-section gives rise to the continuous 
cross-section in ${\underline{\cal B}}_2^{\Delta} $ given by (51),
 as required.
In the case (d), 
 the mapping given by (51) is continuous in the Banach semi-bundle
 ${\underline{\cal B}}_2 ,$ as required.

\begin{flushright}  \large$\diamond$\\
\end{flushright}

\noindent 
\begin{pro} \ 
\begin{description} 
 \item[(1)]   
If $A_p^q$  is constructed under the conditions (a), (b) or (c) of Definition
4.7, then $A_p^q\subseteq L_1 ({\underline {\cal B}}_2^{\Delta } ;
\mu_{H,K}) $. In particular,
if $G/H$ and $G/K$ possess finite invariant measure and  ${1/ p}+{1/{q'}}>1 ,$
then $A_p^q\subseteq L_r ({\underline {\cal B}}_2^{\Delta } ;\mu_{H,K})
 $ where \newline \noindent ${1/ r}={1/ p}+{1/{q'}}-1 .$ \par
\item[(2)] If $A_p^q$  is constructed under the condition (d) and if 
${1/ p}+{1/{q'}}>1 ,$ then \newline \noindent
 $A_p^q\subseteq L_r ({\underline {\cal B}}_2^{\Delta } ;\mu_{H,K})
 $ where  ${1/ r}={1/ p}+{1/{q'}}-1 .$ \par
\item[(3)] If $A_p^q$  is constructed under the condition (d)
 and if ${1/ p}+{1/{q'}}=1 ,$ then 
\newline \noindent $A_p^q\subseteq L_{\infty } ({\underline {\cal B}}_2 ;
\mu_{H\times K}) .$ 

\end{description}
\end{pro}

\noindent
Proof: (1) Consider the space $A_p^q $ under any of the conditions (a) to (c)
given in Definition 4.7. According to the calculations in Proposition 4.6,
we see that 
$$\int_{D(x,y)\in \Gamma }\| \Psi (\sum_{i=1}^{\infty } f_i \otimes g_i )(x,y) \|
d\mu_{H,K} (x,y) \leq \sum_{i=1}^{\infty } \|f_i \|_p \|g_i \|_{q'},$$\noindent for any 
$\sum_{i=1}^{\infty } f_i \otimes g_i \in L_p (\pi )\otimes L_{q' }(\gamma ),$ showing 
that $A_p^q\subseteq L_1 ({\underline {\cal B}}_2^{\Delta} ;\mu_{H,K}) .$

If $G/H$ and $G/K$ have finite invariant measure $\mu_H $ and $\mu_K$
 respectively, we see that  $\mu_{H^x} $ and  $\mu_{K^y}$ are finite invariant
measures on  $G/H^x$ and $G/K^y$ for $x,y\in G .$ Hence 
$\lambda_{H^x } ( z,t) =1=\lambda_{K^y } ( w,t) $ for $z\in G/H^x ,$ 
 $w\in G/K^y$ and $t\in G .$  Using the identity (13), we see that 
$\lambda_{H^x \cap  K^y} =1,$ for almost all $(x,y)\in (G\times G)/
(H\times K). $
\noindent Therefore, for $f\in L_p (\pi ),$ 
$$\int_{G\over {H^x \cap K^y }} \|f(xt)\|^p d\mu_{x,y} (t) = 
\int_{G\over {H^x }}\int_{H^x \over {H^x \cap K^y }}  {\lambda_{H^x \cap  K^y} 
(\alpha,t) \over { \lambda_{H^x} (\alpha ,t)}}\|f(x\alpha t)\|^p 
d\mu_{H_{x,y}} 
(\alpha )d\mu_{H^x } (t), $$
\noindent  where $\mu_{H_{x,y}}$ is the measure on the coset space
 ${H^x / ({H^x \cap K^y })}$ as defined in Lemma 2.1 which is finite 
and invariant (see also the proof of Proposition 4.6). Simplifying, 
\begin{eqnarray*}
\int_{G\over {H^x \cap K^y }} \|f(xt)\|^p d\mu_{x,y} (t) &= &
\int_{G\over {H^x }} \int_{H^x \over {H^x \cap K^y }}\|f(xt)\|^p 
d\mu_{H_{x,y}} (\alpha )d\mu_{H^x } (t) \\
&=& \int_{G\over {H^x }} \|f(xt)\|^p d\mu_{H^x } (t) 
= \|f\|_p^p .
\end{eqnarray*}
\noindent Similarly,
$$\int_{G\over {H^x \cap K^y }} \|g(yt)\|^{q'} d\mu_{x,y} (t) = 
\|g\|_{q'}^{q'} ,$$
\noindent for $g\in L_{q'} .$ Therefore, using Corollary 12.5 of
 Hewitt and Ross 
{\cite{hew}} we obtain
\begin{eqnarray*}
\lefteqn{\int_{G\over {H^x \cap K^y }} \|f(xt)\|\|g(yt)\| d\mu_{x,y} (t)  }
       \\
 & = &
\int_{G\over {H^x \cap K^y }}( \|f(xt)\|^p \|g(yt)\|^{q'} )^{1\over r}
\|f(xt)\|^{1-{p\over r}} \|g(yt)\|^{1-{{q'}\over r}}d\mu_{x,y} (t) , \\
&\leq & \biggl(\int_{G\over {H^x \cap K^y }} \|f(xt)\|^p \|g(yt)\|^{q'}d\mu_{x,y} (t)
\biggr)^{1\over r} \biggl(\int_{G\over {H^x \cap K^y }} \|f(xt)\|^p
 d\mu_{x,y} (t)\biggr)^{{{q'}-1}\over {q'}} \times \\
& & \quad \quad \quad \quad \quad  \quad \quad \quad \quad \quad 
\quad \quad \quad \quad \quad \biggl(\int_{G\over {H^x \cap K^y }}
 \|g(yt)\|^{q'} d\mu_{x,y} (t)\biggr)^{{p-1}\over p} , \\
&=& \biggl(\int_{G\over {H^x \cap K^y }} \|f(xt)\|^p \|g(yt)\|^{q'}d\mu_{x,y} (t)
\biggr)^{1\over r} \|f\|_p^{p({{{q'}-1}\over {q'}})} 
\|g\|_{q'}^{{q'}({{p-1}\over p})} ,
\end{eqnarray*}
which is similar to the right hand side of (42) (Proposition 4.6). 
Note that
\begin{eqnarray*}
\lefteqn{\Biggl(\int_{D(x,y)\in \Upsilon }\| \Psi (\sum_{i=1}^{\infty } f_i 
\otimes g_i )(x,y) \|^r
d\mu_{H,K} (x,y)\Biggr)^{1\over r} }       \\
& \leq & \Biggl(\int_{D(x,y)\in \Upsilon }\biggl(\int_{G\over {H^x \cap K^y }}
\sum_{i=1}^{\infty} \|f_i(xt) \otimes g_i (yt) \| d\mu_{x,y} (t) \biggr)^r
d\mu_{H,K} (D) \Biggr)^{1\over r}\\
& \leq & \Biggl(\int_{D(x,y)\in \Upsilon }\biggl(\sum_{i=1}^{\infty} 
\int_{G\over {H^x \cap K^y }}
 \|f_i(xt)\| \|g_i (yt) \| d\mu_{x,y} (t) \biggr)^r
d\mu_{H,K} (D) \Biggr)^{1\over r}
\end{eqnarray*}
Using the same notation as in (42), by generalised Minkowski's inequality (see Dunford and Schwartz\cite{dun}
p.529) we obtain
\begin{eqnarray*}
\lefteqn{\Biggl(\int_{D(x,y)\in \Upsilon }\| \Psi (\sum_{i=1}^{\infty }
 f_i \otimes g_i )(x,y) \|^r
d\mu_{H,K} (x,y) \Biggr)^{1\over r}}       \\
& \leq & \sum_{i=1}^{\infty } \Biggl(\int_{D(x,y)\in \Upsilon }
(I_i (x,y))^r
d\mu_{H,K} (D) \Biggr)^{1\over r}
\end{eqnarray*}
\noindent Hence using the same calculations which follow inequality (43),
we achieve the required result
 $$\| \Psi (\sum_{i=1}^{\infty } f_i \otimes g_i ) \|_r \leq 
\sum_{i=1}^{\infty}\|f_i\|_p \|g_i\|_{q'}. $$
\par \noindent 
(2) This is evident from (46) (Proposition 4.6).
\par \noindent
(3). Suppose that ${H^x /{(H^x \cap K^y )}}$ and ${K^y / 
{(H^x \cap K^y )}}$ are compact for almost all $(x,y) \in G\times G,$ and 
$p=q.$ Consider the supremum norm on ${\cal B}({\underline {\cal B}}).$ For 
any $\sum_{i=1}^{\infty } f_i \otimes g_i \in L_p (\pi )\otimes L_{p'} (\gamma ^* ),$
\par \noindent 
{\leftline{$\| \Psi (\sum_{i=1}^{\infty }  f_i \otimes g_i ) \|_{\infty } =$}}\par \noindent 
{\rightline{  ess $\sup_{(x,y)\in G\times G}\{ \|\sum_{i=1}^{\infty } 
\int_{G\over {H^x \cap K^y }} {1\over {\lambda (e,t)}}
 \lambda_H (x,t)^{1\over p}f_i(xt)\otimes_{x,y} \lambda_K (y,t)^{1\over {p'}
} g_i(yt)d\mu_{x,y} (t)\|\}. $}}
$\quad $

\noindent Now, following the argument in Proposition 4.6, we see that
$$\| \Psi (\sum_{i=1}^{\infty }  f_i \otimes g_i ) \|_{\infty } \leq {\hbox
{ess }}\sup_{(x,y)}\{M_{x,y}N_{x,y}
  \sum_{i=1}^{\infty }
\|f_i \|_p \|g_i \|_{p'} \} \leq S \sum_{i=1}^{\infty } \|f_i \|_p \|g_i \|_{p'}
,$$
where $S = {\hbox {ess }}\sup_{(x,y)} M_{x,y}N_{x,y} $ is a constant, as required.

\begin{flushright}  \large$\diamond$\\
\end{flushright}

\section{Induced representations, Integral Intertwining Operators and $A_p^q $ spaces}

\subsection{Induced representations and Integral Intertwining Operators}

In this section we shall investigate the possibility of generalising Rieffel's result
(see Rieffel\cite{rief1} Theorem 5.5) on classical $A_p^q $ spaces which 
asserts that such a space is the predual of the space of intertwining 
operators if and only if those operators can be approximated, in the 
ultraweak*-operator topology, by integral operators. To begin, we shall give
 the definition of an integral operator from 
$L_p(\pi) $ to $L_{q} (\gamma),$ and discuss some of its properties. 
\begin{dfn}
 Let $T$ be a bounded linear  operator from $L_p(\pi) $ into $L_{q} (\gamma).$ 
 $T$ is called an \sf integral operator \it if there exists a $\mu_H  \times \mu_K$
measurable function $\Phi ,$ called the \sf kernel \it of $T,$ from $ G/H \times G/K 
$ to ${\cal L}(H(\pi ), H(\gamma ))$
such that for a given $f\in L_p (\pi ),$
\begin{description}
\item[(1)] the function $x\mapsto  \Phi (y,x)f(x)$ is integrable for 
almost all $y\in G/K,$
\item[(2)]  $y\mapsto \int_{G\over H } \Phi (y,x)f(x) d\mu_H (x) $ 
belongs to $L_q (\gamma )$ and 
\item[(3)] $(Tf)(y) = \int_{G\over H}  \Phi (y,x)
f(x)d\mu_H (x),$ for almost all $y\in G/K. $
\end{description}
\end{dfn}

The next result describes the properties of the kernel of an intertwining 
integral operator. The existence of such operators will be discussed in
Proposition 5.3.

\begin{pro}
Let $\Phi $ be the kernel of a given integral intertwining operator for induced representations
$U_p^{\pi }$ and $U_{q}^{\gamma } .$  Then $\Phi $ satisfies the
following properties.
\begin{description} 
\item[(1)] For almost all  $x\in G/H ,$ $y\in G/K$ and for all $ s\in G ,$ 
\begin{eqnarray}
\lambda_H (x,s^{-1})^{1\over p'}\Phi (y,xs^{-1}) &= &
\lambda_K (y,s)^{1\over q}
\Phi (ys,x). 
\end{eqnarray}
\item[(2)] For all  $h\in H,$ $k\in K,$ and for almost all $x\in G/H,$
 $y\in G/K,$ 
\begin{eqnarray}
\Phi (ky,hx){\pi}_h &=&  {{\gamma}_k}\Phi (y,x).
\end{eqnarray}   
\item[(3)] Under the conditions given in  
Definition 4.7, $\Phi (y,x)$ is an intertwining operator of the
representations 
 ${\pi}^x$ and ${\gamma}^y$ of the subgroup $H^x \cap K^y$ of $G$
for  almost all $x\in G/H $ and  $y\in G/K.$
\end{description}
\end{pro}

\par \noindent                                                        
Proof : (1) Suppose that $T$ is an integral operator from $ L_p (\pi )$ to
$ L_q (\gamma )$ with the kernel $\Phi .$  Then for $f\in L_p (\pi )$ and $y\in G,$
 $$(Tf)(y) = \int_{G\over H} \Phi (y,x) f(x) d\mu_H (x).$$
 In addition, if
$T\in Hom_G (L_p (\pi),L_{q} (\gamma))$ then
$$(TU^{\pi}_s f)(y) = (U^{\gamma}_s T f)(y){\hbox { for almost all }}
y\in G/K{\hbox { and for }}s\in G.$$
Now 
\begin{eqnarray*}
(TU^{\pi}_s f)(y) &=& \int_{G\over H} \Phi (y,x) \lambda_H (x,s)^{1\over p}
 f(xs) d\mu_H (x).
\end{eqnarray*}

\noindent Changing variables $xs \mapsto x$, we find 
$$(TU^{\pi}_s f)(y) =\int_{G\over H} \Phi (y,xs^{-1})\lambda_H 
(x,s^{-1})\lambda_H (xs^{-1},s)^{1\over p} 
f(x)d\mu_H (x).$$
Since $\lambda_H (x,s^{-1})\lambda_H (xs^{-1},s)=1,$ the
 above integral 
simplifies  to
\begin{eqnarray}
(TU^{\pi}_s f)(y)&=&\int_{G\over H} \Phi (y,xs^{-1})\lambda_H (x,s^{-1})^{1\over p'} 
f(x)d\mu_H (x).
\end{eqnarray}
On the other hand,
\begin{eqnarray}
(U^{\gamma}_s Tf)(y)
 &= &\lambda_K(y,s)^{1\over q}\int_{G\over H} \Phi (ys,x)f(x)d\mu_H (x) . 
\end{eqnarray}
\noindent
Therefore, by (55) and (56), property (1) follows.

\noindent
(2) For $k\in K$ and $y\in G,$
\begin{eqnarray}
\gamma _k (Tf)(y) = (Tf)(ky) &= &\int _{G\over H} \Phi (ky,hx)\pi _h f(x)d\mu _H (x), 
\end{eqnarray}
for  $h\in H .$ On the other hand,\\
\begin{eqnarray}
\gamma _k (Tf)(y)& =& \gamma _k \int _{G\over H} \Phi (y,x) 
f(x)d\mu _H (x) .
\end{eqnarray}
It is clear that property (2) follows from  (57) and (58).

\noindent
(3) We want to show that
\begin{eqnarray}
\quad \quad \quad \quad {\gamma_b^y}\Phi (y,x)&=& \Phi (y,x) {\pi_b^x}  ,
\end{eqnarray}
for all $b\in H^x \cap K^y$  and for almost all $x\in G/H $ and $y\in G/K.$
For any $b\in H^x \cap K^y$ we have
 $b=y^{-1}ky =x^{-1}hx $ for some $ h\in H $ and $ k\in K.$
Using (54),
\begin{eqnarray*}
{\gamma_{yby^{-1}}} \Phi (y,x)&=&\Phi ({yby^{-1}}y,{xbx^{-1}}x) 
{\pi_{xbx^{-1}}} ,
\end{eqnarray*}
which implies
\begin{eqnarray}
{\gamma_b^y}  \Phi (y,x)&=&  \Phi (yb,xb) {\pi_b^x},
= {{\lambda_H (xb,b^{-1} )^{1\over {p'}}}\over 
{\lambda_K (y,b )^{1\over {q}}} } \Phi (y,x) {\pi_b^x}, {\hbox { by
(48) ,}} \nonumber \\
&=& {1\over {{\lambda_H (x,b )^{1\over {p'}}}
{\lambda_K (y,b )^{1\over q}}} } \Phi (y,x) {\pi_b^x}, 
\end{eqnarray}
Under  conditions (a), (b) or (c) of Def.4.7, (60) simplifies to
(59), as required.

Now suppose that the condition given in (d)  of Definition 4.7
applies. Consider the right hand side of (60). We see that
\begin{eqnarray*}
 {1\over {{\lambda_H (x,b )^{1\over {p'}}}
{\lambda_K (y,b )^{1\over q}}} } &=& \biggl( {{\lambda_{H^x \cap K^y} (e,b)}
\over {\lambda_K (y,b )}}\biggl)^{1\over {p'}}\biggl( {{\lambda_{H^x \cap K^y} (e,b)}
\over {\lambda_H (x,b )}}\biggl)^{1\over q}, {\hbox{ by (9),}} \\
&=& \biggl( {{\lambda_{H^x \cap K^y} (e,b)}
\over {\lambda_{K^y } (e,b )}}\biggl)^{1\over {p'}}\biggl( {{\lambda_{H^x \cap K^y} (e,b)}
\over {\lambda_{H^x} (e,b )}}\biggl)^{1\over q}, 
\end{eqnarray*} 
Under the condition that $H^x /(H^x \cap K^y )$ and $K^y/(H^x \cap K^y )$
have  invariant measure, we have 
$${{\lambda_{H^x \cap K^y} (e,b)}
\over {\lambda_{H^x } (e,b )}} = {{\lambda_{H^x \cap K^y} (e,b)}
\over {\lambda_{K^y } (e,b )}} =1 $$
(see (32) in the proof of Proposition 4.5). Therefore,
\begin{eqnarray}
{\gamma_b^y}  \Phi (y,x)&=&  \Phi (y,x) {\pi_b^x},
\end{eqnarray} 
for all $b\in H^x \cap K^y$  and for almost all $x\in G/H $ and $y\in G/K.$
Hence the result.

\begin{flushright}  \large$\diamond$\\
\end{flushright}

Following an argument similar to that of Moore\cite{moo}, we shall obtain a 
result for intertwining operators between $L_1 (\pi )$ and
 $L_q (\gamma ),q>1.$

\par \noindent
\begin{pro}
 Let  $U_1^{\pi }$ and $U_q^{\gamma }$
 be induced representations of the locally compact group $G$ with the 
corresponding Banach spaces of functions $L_1 (\pi )$ and $L_q (\gamma )$
\newline \noindent $ (q>1)$, respectively. Then,  if the Banach space $\cal {H}(\pi )$ is
 separable, the intertwining operators $T$ for these 
representations  are integral operators with the corresponding  kernel 
$\Phi $ satisfying
\par $ess$ $sup_{x\in {G\over H} } \bigl (\int_{G\over K } \| \Phi (y,x )\|^q d\mu_K (y)\bigr )^{
1\over q} \leq \|T \|.$
\end{pro}

 \noindent Proof : The proof is in two parts:
\par \noindent (1). Let $S$ and $R$ be  fixed Borel cross-sections
of $H$ and $K$ in $G$. Then $G/H \simeq S,$  $G/K \simeq R$ and 
we regard $\mu_H $ and $\mu_K $ as measures on $S$ and $R$ .
Let $C$ be a continuous linear map of  $L_1 (S,\cal {H}(\pi ), \mu_H )$ 
into 
$L_q (R,\cal {H}(\gamma ),\mu_K )$. Firstly we prove that $C$ is an integral
operator.
 For $u\in H(\pi )$ define
$C_u : L_1 (S,\mu_H ) \rightarrow L_q (R,\cal {H}(\gamma ),\mu_K ) $ by
\begin{eqnarray}
\quad \quad \quad \quad \quad \quad \quad \quad C_u (g) &= & C(gu),
\end{eqnarray}
 for $ g\in L_1 (S,\mu_H ) .$ $C_u (g)$ is bounded since 
 $\| C_u \| \leq \| C \|.\| u \|$.
Then by Dunford  and Schwartz\cite{dun}, Theorem 10, p.507, there 
exists a $\mu_H $-essentially
unique bounded measurable function $\chi_u $(.) on $S$ to a weakly
compact subset of $L_q (R,\cal {H}(\gamma ),\mu_K )$ such that 
$$C_u (g) = \int_S \chi_u (s)g(s) d\mu_H (s)  
$$ and $\| C_u \| =$ ess $\sup \| \chi_u (s) \|$.
Let $K_u (.,s) = \chi_u (s)$ so that $K_u :R\times S \mapsto H(\gamma).$
Then $K_u $ is $\mu_H \times \mu_K $ measurable (see Dunford  and 
Schwartz\cite{dun}, Theorem 17 p.198), and we have
\begin{eqnarray}
\quad \quad \quad \quad (C_u (g))(t) &=& \int_S g(s)K_u (t,s) d\mu (s)
\end{eqnarray}
 with ess $\sup_{s\in S } \bigl (\int_R   \| K_u (t,s )\|^q d\mu_K (t)\bigr
 )^{1\over q} \leq \| C\|\|u\|.$ 
Following the same argument as in Moore\cite{moo} we can define a map $K$ 
on $R\times S$  in to the space of bounded linear operators from 
${\cal {H}}(\pi )$ to ${\cal {H}}(\gamma )$ such that $K (t,s) =0$ for 
$(t,s)  $ in a suitably chosen null set $N$
 and $K (t,s)u
 = K_u (t,s)$ otherwise, for each $u\in {\cal {H}}(\pi ),$ with $\|K(s,t)\|\leq C.$ Then, by
 (62) and (63), we have
$$C(gu)(t) = \int_S K (t,s)g(s)u d\mu (s)$$and 
ess $\sup_{s\in S } \bigl (\int_R  {{\| K(t,s )u\|}\over {\|u\|}}^q 
d\mu_K (t)\bigr )^{1\over q} \leq \| C\|{\hbox { for any }} u\in H(\pi ),$  
which implies that ess $\sup_{s\in S } \bigl (\int_R  \| K (t,s )\|^q d\mu_K (t)
\bigr )^{1\over q} \leq \| C\|.$
Hence for $g\in L_1 (S,\cal {H}(\pi),\mu_H) $ we have
\begin{eqnarray}
(Cg)(t) = \int_S K(t,s)g(s) d\mu_H (s).
\end{eqnarray}

\par \noindent
 (2). Secondly, we prove that the intertwining operators $T$ from
 $L_1 (\pi)$ to $L_q (\gamma)$ are integral 
operators.

\par  Observing that $G\simeq H\times S,$  for a given continuous function
$f' \in L_1 (S,{\cal {H}}(\pi ),\mu_H)$ we can define a function 
$ (\Phi_1 {f'})\in L_1 (\pi )$ by
$$ (\Phi_1 {f'})(y) = \pi (h)f'(s),$$ 
where  $y\in G$ with  $ y=hs$ for  $ h\in H$ 
and $s\in S.$  Similarly, since $G\simeq K\times R,$ for a given continuous
function $g\in L_q (\gamma) ,$ we define the function $(\Phi_q g )\in 
L_q (R,\cal {H}(\gamma ), \mu_K)$ by 
$$(\Phi_q g ) (r) = g(r)$$ 
for $r\in R .$ Clearly, $\|f' \|_1 = \|(\Phi_1 {f'})\|_1$ and $\|g \|_q =
\|(\Phi_q g )\|.$

 For a given intertwining 
operator $T$ from  $L_1 (\pi)$ to $L_q (\gamma)$ we define an operator
${\tilde{T}} $ on the space of continuous functions in 
$L_1 (S,{\cal {H}}(\pi ),\mu_H)$ to $ L_q (R,
{\cal {H}}(\gamma ),\mu_K)$ by
$${\tilde{T}} := \Phi_q T \Phi_1 .$$ 
Since the space of continuous functions in $L_1 (S,{\cal {H}}(\pi ),\mu_H)$
is dense in $L_1 (S,{\cal {H}}(\pi ),\mu_H),$  
we have the following 
commutative diagram:
$$L_1 (S,{\cal {H}}(\pi ),\mu_H) \stackrel{{\tilde{T}}}{\mapsto } L_q (R,
{\cal {H}}(\gamma ),\mu_K)$$
$$\Phi_1 \downarrow \qquad \quad \uparrow \Phi_q$$
$$L_1 (\pi) \stackrel{T}{\mapsto } L_q (\gamma)$$ 
with ${\tilde{T}} (f') = \Phi_q T \Phi_1 (f')$ for $f' \in 
L_1 (S,{\cal {H}}(\pi ),\mu_H) .$ Clearly, $\|T\|= \|{\tilde{T}}\|.$ 
Using the result in part (1), we see that there  exists a map
  $K$ from $S\times R $ to the set of bounded linear maps from ${\cal H }
(\pi )$ to ${\cal H } (\gamma )$ such that $$({\tilde{T}}f')(t) = \int_S K(t,s) f'(s) 
d\mu_H (s), $$ for $f' \in L^1 (S,\cal {H}(\pi ),\mu_H)$ and $t\in R.$
Using the Borel isomorphism $G\simeq K\times R$ any $y\in G$ can be written as
$y= k(e,y) \ell (e,y)$ where $k(e,y) \in K$ and $\ell (e,y) \in R.$ Both $k$ and $\ell $
are Borel functions on $R\times G.$
Then, for $f\in L_1 (\pi) $ 
\begin{eqnarray*}
(Tf)(y) = ({\Phi_q}^{-1} {\tilde{T}} {\Phi_1}^{-1} (f))(y)& = &
\gamma (k(e,y)) (({\tilde{T}} {\Phi_1}^{-1})f)(\ell (e,y)), \\
  & = &\gamma (k(e,y)) \int_S K(\ell (e,y),s)(( {\Phi_1}^{-1})f)(s) d\mu_H (s).
\end{eqnarray*}
 But since $({\Phi_1 }^{-1} f)(s)= f(s)$, we have
\begin{eqnarray*}
\quad \quad \quad \quad (Tf)(y) &= &\gamma (k(e,y)) \int_S K(\ell (e,y),s) f(s) d\mu_H (s).
\end{eqnarray*}
Now the Borel isomorphism $G\simeq H\times S,$ allows us to express any $x\in G$ in the form
$x= h(e,x) m(e,x),$ where $h$ and $m$ are Borel functions on $H\times S,$  $h(e,x)\in H$ 
and $m(e,x)\in S.$
If we define $\Phi (y,x) = \gamma (k(e,y))K(\ell (e,y),m(e,x)){\pi (h(e,x)) }^* ,$ then
 we have $\|\Phi (y,x)\| = \| K(\ell (e,y),m(e,x))\|$ and 
\begin{eqnarray}
(Tf)(y) = \int_S \Phi (y,s) f(s) ds & = &\int_{G\over H} 
\Phi (y,x) f(x) d{\mu_H }(x),
\end{eqnarray}
 with ess $\sup_{x\in G/H } \bigl (\int_{G/K } \| \Phi (y,x )\|^q 
d\mu_K (y)\bigr )^{1\over q} \leq \|T \|.$
Therefore T is an integral operator.

\begin{flushright}  \large$\diamond$\\
\end{flushright}
\subsection{The space $A_p^q $ as the predual of the space of 
intertwining operators}

We are now in a position to state the main result of this section, which is a 
generalisation of Rieffel's result(\cite{rief1} Theorem 5.5) on classical $
A_p^q $ spaces.

 \noindent 
\begin{thm}
Suppose that the space $A_p^q ,$ $({q'}>1,)$ is constructed under one of the conditions given
in Definition 4.7. Then the following statements are equivalent. 
\begin{description}
\item[(a)] $L_p (\pi) \otimes_G^{\sigma }  L_{q'} (\gamma^*)  \simeq  A_p^q .$ 
\item[(b)]  Every element of $Int_G (U_p^{\pi } , U_q^{\gamma })$ can be approximated 
in the ultraweak*-operator topology by integral operators.
\end{description}
\end{thm}

\noindent 
Proof:(b)$\Rightarrow $(a)
 Suppose that every element of  $Int_G (U_p^{\pi } , U_q^{\gamma })$ can be approximated 
in the ultraweak*-operator topology by integral operators.
 First we show that the kernel of $\Psi$ contains the subspace  $L$  of
 $L_p (\pi) \otimes  L_{q'} (\gamma^*) ;$ 
  that is,
$$\Psi(\Sigma_{i=1}^{\infty}U^{\pi} (s)f_i \otimes g_i) = \Psi(\Sigma_{i=1}^
{\infty}f_i \otimes {(U^{\gamma})^*} (s) g_i)$$  for $s\in G$. 
 In the following we write $\lambda (\cdot ,\cdot )$ for 
$\lambda_{H^x\cap K^y} (\cdot ,\cdot ).$
\noindent Now 
\begin{eqnarray*}
\lefteqn{\Psi(\sum_{i=1}^{\infty}U^{\pi}(s)f_i \otimes g_i)(x,y) }     \\
 &= &\int_{G\over {H^x\cap K^y} }\sum_{i=1}^{\infty} {1\over { \lambda  (e,t)} } 
\lambda_H (x,t)^{1\over p}\lambda_H (xt,s)^{1\over p} f_{i}(xts)
\otimes_{x,y}  \lambda_K (y,t)^{1\over {q'}}g_{i}(yt) d\mu_{x,y}(t),\\ 
&= &\int_{G\over {H^x\cap K^y}}\sum_{i=1}^{\infty} {1\over { \lambda  (e,t)} } 
\lambda_H (x,ts)^{1\over p} f_{i}(xts)\otimes_{x,y}  
\lambda_K (y,t)^{1\over {q'}}g_{i}(yt) d\mu_{x,y}(t), \\     
& =&\int_{G\over{H^x\cap K^y}}\sum_{i=1}^{\infty} {\lambda (t,s^{-1})\over 
{\lambda (e,ts^{-1})}}\lambda_H (x,t)^{1\over p} f_{i}(xt)\otimes_{x,y} 
\lambda_K (y,ts^{-1})^{1\over {q'}}g_{i}(yts^{-1}) d\mu_{x,y}(t),
\end{eqnarray*}
\noindent on changing variables $ts\mapsto t.$
\noindent Since ${\lambda (t,s^{-1})/ {\lambda (e,ts^{-1})}}= {1/ {
\lambda (e,t)}},$  and \newline \noindent 
$\lambda_K (y,ts^{-1}) = \lambda_K (yt,s^{-1})
\lambda_K (y,t),$
\noindent
\begin{eqnarray*}
\lefteqn{\Psi(\sum_{i=1}^{\infty}U^{\pi}(s)f_i \otimes g_i)(x,y) }  \\
 &=&\int_{G\over {H^x\cap K^y}}\sum_{i=1}^{\infty} {1\over { \lambda  (e,t)}
 } \lambda_H (x,t)^{1\over p}
 f_{i}(xts)\otimes_{x,y} \\
& &\quad \quad \quad \quad \quad \quad \quad \quad \quad \quad \quad \quad 
 \lambda_K (y,t)^{1\over {q'}}
\lambda_K (yt,s^{-1})^{1\over {q'}}g_{i}(yts^{-1})
  d\mu_{x,y}(t),\\
 &=&\int_{G\over{H^x\cap K^y}} \sum_{i=1}^{\infty}{1\over { \lambda  (e,t)} }
 \lambda_H (x,t)^{1\over p}
  f_i(xt) \otimes_{x,y} (U^{{\gamma}})^*(s) g_i(yt)d\mu_{x,y} (t),\\
& =&\Psi (\sum_{i=1}^{\infty} f_i\otimes (U^{\gamma})^* (s) g_i).
\end{eqnarray*}

Now it only requires to prove that the kernel of $\Psi$ is contained in $ 
L.$
To achieve this, it suffices to
show that any bounded linear functional $F$ on $L_p (\pi) 
\otimes_G^{\sigma } L_q (\gamma)$
which annihilates $L$ also annihilates the kernel of $\Psi$. Since $F$ annihilates 
$L,$ there exists $T\in Int_G (U_p^{\pi } , U_q^{\gamma })$ such that
\begin{eqnarray}
\quad \quad \quad \quad \quad \quad \langle r,F\rangle  &=& \sum_{i=1}^{\infty} \langle g_i,Tf_i\rangle ,
\end{eqnarray}
 for any $r\in L_p (\pi) \otimes_G^{\sigma } L_{q'} (\gamma^*)$ with the expansion 
$$r = \sum_{i=1}^{\infty} f_i\otimes g_i .$$

\noindent
Suppose now that $r$ is in the kernel of $\Psi$. Then,
\begin{eqnarray}
\sum_{i=1}^{\infty}\int_{G\over {H^x\cap K^y}}{1\over 
{ \lambda  (e,t)} } \lambda_H (x,t)^{1\over p} f_i (xt)\otimes_{x,y}
\lambda_K (y,t)^{1\over {q'}} g_i (yt)\mu_{x,y}(t) =0.
\end{eqnarray}
 By (66), it suffices to show that  
$$\sum_{i=1}^{\infty} \langle g_i,Tf_i\rangle =0.$$
 Under the assumption that the operator $T$ can be approximated by the
integral operators $\{T_j : j\in I\}$
 in the ultraweak*-operator topology, we have
 $$\sum_{i=1}^{\infty } \langle g_i , T_j f_i \rangle {\rightarrow } \sum_{i=1}
^{\infty } \langle g_i , T f_i \rangle .$$
 Hence in order to prove $\sum_{i=1}^{\infty} \langle g_i,Tf_i\rangle =0,$ it is 
sufficient to prove
 $$\sum_{i=1}^{\infty } \langle g_i , T_j f_i \rangle =0, $$
\noindent for each $T_j .$  Since $  T_j $ is an integral operator, we have
$$(T_j f_i)(y) = \int_{G\over H} \Phi_j (y,x)f_i(x) d\mu_H (x) ,$$ where 
$ \Phi_j $ is the kernel of $T_j$ as described in  Definition 5.1.
Thus,
\begin{eqnarray*}
\lefteqn{\sum_{i=1}^{\infty} \langle g_i,T_j f_i\rangle }\\
 &=& \sum_{i=1}^{\infty}\int_{G\over K} \langle g_i(y),(T_j f_i)(y)\rangle
 d\mu_K (y),
 \\
 &=& \sum_{i=1}^{\infty}\int_{G\over K} \int_{G\over H} \langle g_i(y),
\Phi_j (y,x)f_i(x)\rangle d\mu_H (x) d\mu_K (y), \\
 &=&\sum_{i=1}^{\infty} \int_{{G\times G}\over {H\times K}} \langle g_i (y),
\Phi_j (y,x)f_i (x)\rangle d\mu_{H\times K} (x,y), \\
 &=&\sum_{i=1}^{\infty} \int_{D\in \Upsilon } \int_{G \over 
{H^x \cap K^y}}\langle g_i (yt),\Phi_j (yt,xt)f_i(xt)\rangle
 d\mu_{x,y}
(t)d\mu_{(H,K)}(D),  
\end{eqnarray*}
\noindent using disintegration of measures as explained in Lemma 2.2. 
(Also, see the discussion preceeding the Lemma).
By Proposition 5.2 (1), $\lambda_H (xt,t^{-1})^{1\over p'}\Phi_j (y,x) = \lambda_K 
(y,t)^{1\over q}\Phi_j (yt,xt)$  for almost all $x\in G/H.$ 

 \noindent Therefore, \par \noindent 
\begin{eqnarray*}
\lefteqn{\sum_{i=1}^{\infty} \langle g_i,T_j f_i\rangle } {\hspace{0.1in}}\\
&=&\sum_{i=1}^{\infty}\int_{D\in \Upsilon } 
\int_{G \over {H^x \cap K^y}}\langle g_i (yt),
{\lambda_H (xt,t^{-1})^{1\over p'}
\over {\lambda_K(y,t)^{1\over q}}}\Phi_j (y,x)f_i (xt)\rangle 
 d\mu_{x,y}(t) d\mu_{(H,K)}(D) .
\end{eqnarray*}
\noindent From the identity (13), we see that 
$${\lambda_H (xt,t^{-1})^{1\over p'}\over {\lambda_K(y,t)^{1\over q}}}=
{1\over{\lambda_H (x,t)^{1\over p'} \lambda_K(y,t)^{1\over q}}}=
{1\over {\lambda (e,t)}}
\lambda_H (x,t)^{1\over p} \lambda_K(y,t)^{1\over q'}.$$  
\noindent Consequently, 
\begin{eqnarray}
\lefteqn{\sum_{i=1}^{\infty} \langle g_i,T_j f_i\rangle }\nonumber \\ 
&= & \sum_{i=1}^{\infty} \int_{D\in \Upsilon } 
\int_{G \over {H^x \cap K^y}}{1\over{\lambda (e,t)}}
\langle \lambda_K(y,t)^{1\over q'} g_i (yt),\Phi_j (y,x)\lambda_H (x,t)^{1\over p}
 f_i (xt)\rangle \nonumber \\ 
& & \quad \quad \quad \quad \quad \quad \quad \quad \quad 
\quad \quad \quad \quad \quad \quad \quad d\mu_{x,y}(t)d\mu_{(H,K)} (D).
\end{eqnarray}

\noindent By Proposition 5.2 (3), $\Phi_j (y,x)\in Int_{H^x \cap K^y}
(H(\pi^x),H(\gamma^y))$  under the conditions given in Definition
4.7.  Hence
there exists \newline \noindent $ \Theta_j (y,x) \in 
(H(\pi^x)\otimes_{H^x \cap K^y} H(({\gamma^y})^*))^* $  such that 
\begin{eqnarray*}
\lefteqn{\sum_{i=1}^{\infty}\langle \lambda_K(y,t)^{1\over q'} g_i (yt),\Phi_j 
(y,x)\lambda_H (x,t)^{1\over p} f_i(xt)\rangle
 }{\hspace{0.3in}} \\
&=&\sum_{i=1}^{\infty}\langle \lambda_H (x,t)^{1\over p} 
f_i (xt) \otimes_{x,y}\lambda_K(y,t)^{1\over q'}
  g_i (yt),\Theta_j (x,y)\rangle ,
\end{eqnarray*}
\noindent (see Rieffel\cite{rief2}). Therefore we have,
\begin{eqnarray}
\lefteqn{\sum_{i=1}^{\infty }\langle g_i,T_j f_i \rangle = }\nonumber \\
& &\sum_{i=1}^{\infty} \int_{D\in  \Upsilon }\int_{\Delta 
\over{(H\times K)^{(x,y)}\cap \Delta}}
\langle {1\over{\lambda (e,t)}}\lambda_H (x,t)^{1\over p}f_i (xt)\otimes_{x,y}
\lambda_K(y,t)^{1\over q'}g_i (yt),\Theta_j (x,y)\rangle \nonumber \\
& & \quad \quad \quad \quad \quad \quad \quad \quad \quad \quad 
 \quad \quad \quad \quad \quad \quad \quad \quad d\mu_{x,y}(t) 
 d\mu_{(H,K)}(D). 
\end{eqnarray}
 \noindent 
 Hence, by (67), $$\sum_{i=1}^{\infty }\langle g_i ,T_j f_i\rangle = 0 ,$$
 as required.
 \par \noindent (a)$\Rightarrow $ (b)
 Now suppose that the kernel of $\Psi $ is $L.$ We want to show that the integral 
 operators of the form $T_{\phi } f(y) = \int_{G/H}  \phi (y,x) f(x) d\mu_H (x) $ form a dense set 
 in $Hom_G (L_p (\pi ) , L_q (\gamma ))$ in the ultraweak*-operator topology; or equivalently, 
the corresponding linear functionals are dense in $(L_p (\pi ) \otimes_G L_{q'} (\gamma^* ))^* $
in the weak*- topology. Hence, we only need to show that the annihilator of these functionals, 
regarded as functionals on  $(L_p (\pi ) \otimes_{\sigma } L_{q'} (\gamma^* ))^* , $ is $L.$  
But by (69) we see that the annihilator  of these linear functional is the kernel of $\Psi$ which is 
equal to $L$ under our assumption. This concludes the proof of the Theorem.

\begin{flushright}  \large$\diamond$\\
\end{flushright}

\begin{cor} Suppose that every element of $Int_G (U_p^{\pi } , U_q^{\gamma })$
 can be approximated 
in the ultraweak*-operator topology by integral operators.
Then the intertwining number  $\partial (U_p^{\pi },U_q^{\gamma })$ is 
equal to the dimension of the space of all
 functions $\Phi $ given in  Definition 5.1. Moreover, if $H$ and $K$ are
discretely related, 
$$\partial (U_p^{\pi },U_q^{\gamma }) = \sum_{\vartheta \in \Upsilon } d_{\vartheta
} ,$$ where $d_{\vartheta } $ is the dimension of the set of all functions 
$\Phi $ which vanish outside the double coset ${\vartheta }$.
\end{cor}

\noindent
Proof: Let $T\in Int_G (U_p^{\pi } , U_q^{\gamma }).$ By (69) we have
\begin{eqnarray}
\sum_{i=1}^{\infty }\langle g_i,Tf_i \rangle & =& \int_{D\in \Upsilon} \langle \Psi (x,y) , \Theta (x,y)
\rangle d\mu_{(H,K)}(D)  = \langle \Psi , \Theta \rangle .
\end{eqnarray}

\noindent Now using Proposition 1.5 and Theorem 5.4, 
$$(A^q_p)^* \simeq Hom_G (L_p (\pi ),{L}_q (\gamma )).$$
By (70), 
the intertwining number  $\partial (U_p^{\pi },U_q^{\gamma })$ is 
equal to the dimension  of the space 
of all functions $\Theta $ which, in turn is equal to  the dimension  of the space 
of all functions $\Phi. $ 
\par 
If $H$ and $K$ are discretely related, $G$ is a union of 
a null set and a countable collection of double cosets. By Proposition 5.2
 (2), the value  of $\Phi $  on $\vartheta $ is uniquely determined
 by its value $\Phi (x_0, y_0)$
 at $(x_0, y_0)$ where $(x_0, y_0) \in \vartheta $.
\par \noindent Hence $$\partial (U_p^{\pi },U_q^{\gamma }) = \sum_{\vartheta \in D} d_{\vartheta }.$$

\begin{flushright}  \large$\diamond$\\
\end{flushright}

\end{document}